\documentclass[sn-mathphys,Numbered]{sn-jnl}% Math and Physical Sciences Reference Style
\usepackage{graphicx}%
\usepackage{multirow}%
\usepackage{amsmath,amssymb,amsfonts}%
\usepackage{amsthm}%
\usepackage{mathrsfs}%
\usepackage[title]{appendix}%
\usepackage{xcolor}%
\usepackage{textcomp}%
\usepackage{manyfoot}%
\usepackage{booktabs}%
\usepackage{algorithm}%
\usepackage{algorithmicx}%
\usepackage{algpseudocode}%
\usepackage{listings}%
\def\lvtimes{\vec{\ltimes}}

\def\J{{\bf 1}}
\DeclareMathOperator{\sgn}{sgn}

\DeclareMathOperator{\Col}{Col}

\DeclareMathOperator{\lcm}{lcm}

\DeclareMathOperator{\gl}{gl}

\def\cal{\mathcal}

\def\com{com}
\def\ra{\rightarrow}

\def\a{\alpha}
\def\b{\beta}
\def\d{\delta}

\def\D{\Delta}

\def\0{{\bf 0}}

\def\ID{\mathrm{ID}}
\newcommand{\R}{{\mathbb R}}
\newcommand{\C}{{\mathbb C}}
\newcommand{\N}{{\mathbb N}}

\newcommand{\F}{{\mathbb F}}

\def\dsum{\mathop{\sum}\limits}

\theoremstyle{thmstyleone}%
\newtheorem{thm}{Theorem}%  meant for continuous numbers
\newtheorem{prp}[thm]{Proposition}% 

\newtheorem{cor}[thm]{Corollary}

\theoremstyle{thmstyletwo}%
\newtheorem{exa}{Example}%
\newtheorem{rem}{Remark}%

\theoremstyle{thmstylethree}%
\newtheorem{dfn}{Definition}%

\begin{document}
	
	\title[Contracted Product of Hypermatrices via STP of Matrices]{Contracted Product of Hypermatrices via STP of Matrices}
	
	
	\author[1]{\fnm{Daizhan} \sur{Cheng}}\email{dcheng@iss.ac.cn}
	
	\author[2]{\fnm{Min} \sur{Meng}}\email{mengmin@tongji.edu.cn}
%	\equalcont{These authors contributed equally to this work.}
	
	\author[1,3]{\fnm{Xiao} \sur{Zhang}}\email{xiaozhang@amss.ac.cn}
%	\equalcont{These authors contributed equally to this work.}
	
	\author*[1,4]{\fnm{Zhengping} \sur{Ji}\email{jizhengping@amss.ac.cn}}
%	\equalcont{These authors contributed equally to this work.}
	
	\affil[1]{\orgdiv{Key Laboratory of Systems and Control}, \orgname{Academy of Mathematics and Systems Science, Chinese Academy of Sciences}, \orgaddress{\city{Beijing}, \postcode{100190}, \country{P. R. China}}}
	
	\affil[2]{\orgdiv{Department of Control Science and Engineering, School of Electronics and Information Engineering}, \orgname{Tongji University}, \orgaddress{\city{Shanghai}, \postcode{201210}, \country{P. R. China}}}
	
	\affil[3]{\orgdiv{National Center for Mathematics and Interdisciplinary Sciences}, \orgname{Chinese Academy of Sciences}, \orgaddress{\city{Beijing}, \postcode{100190}, \country{P. R. China}}}
	
	\affil[4]{\orgdiv{School of Mathematical Sciences}, \orgname{University of Chinese Academy of Sciences}, \orgaddress{\city{Beijing}, \postcode{100149}, \country{P. R. China}}}
	
	
	\abstract{An equivalent definition of hypermatrices is introduced. The matrix expression of hypermatrices is proposed. Using permutation matrices, the conversion of different matrix expressions is revealed. The various contracted products of hypermatrices are realized by semi-tensor products (STP) of matrices via matrix expressions of hypermatrices.}
	
		\footnotetext{This work is supported partly by the National Natural Science Foundation of China (NSFC) under Grant 62073315 and 62103305, and China Postdoctoral Science Foundation 2021M703423 and 2022T150686.}
	
	
	
	\maketitle
	
	\section{Introduction}\label{sec1}
	
	The hypermatrix is an extension of the matrix to the higher order case (order $d\geq 3$). It has found wide application in many fields, including computer science \cite{vas03}, signal processing \cite{del01}, statistics \cite{del00}, etc. We refer to \cite{lim13} for a systematic introduction to the hypermatrices, and to \cite{bus19,fan20} for some later developments.
	
	A generalized matrix product, called the semi-tensor product (STP), was proposed two decades ago \cite{che01}.  Since then it has been rapidly developed both in theoretical aspects and in various applications \cite{che11,che12}. For example, it has been applied to the study of Boolean and finite-valued networks (see survey papers \cite{for16,li18,lu17,muh16}); finite games (see survey paper \cite{che21}); finite automata (see survey paper \cite{yan22}); dimension-varying systems \cite{che19,che19b}, etc.
	
	Later, in addition to the original matrix-matrix STP, the matrix-vector STP was proposed \cite{che19c} and applied to dimension-varying dynamic (control) systems \cite{che19b}. Recently, the STP of hypermatrices has also been introduced \cite{che23}. It seems possible that some new unknown STPs will appear in the future. In the early days, the (matrix-matrix) STP is said to be a generalization of the conventional matrix product. With the appearance of new STPs, this explanation does not seem to be sufficient.
	
	The purpose of this article is to show that in essence, STPs are representations of multilinear mappings over hypermatrices. More precisely, the hypermatrices are first expressed in their matrix forms. Then the multilinear mappings over them can be executed by STP over their matrix expressions. For this purpose, the conversion of different matrix expressions of hypermatrices plays an important role.
	
	The rest of this paper is organized as follows: Section 2 considers how to transform a hypermatrix into its distinct matrix expressions. The permutation matrix is also introduced for transforming among different matrix expressions. Section 3 shows how to realize the contracted product of hypermatrices over their matrix expressions. In Section 4, the STP is used to perform the contraction product of hypermatrices, which shows that STPs are essentially multilinear operators over hypermatrices.
	
	Before ending this section, we give a list of the notations used in this paper.
	
	\begin{enumerate}
		\item $\N$: the set of positive integers.
		\item $\R$: the set of real numbers.
		\item $\C$: the set of complex numbers.
		\item $\F^n$: the $n$ dimensional Euclidean space over $\F$, where $\F$ is a field (Particularly, $\F=\R$ or $\F=\C$).
		\item $\lcm(a,b)$: the least common multiple of $a$ and $b$, ($a,b\in \N$).
		\item  $\F^{n_1\times n_2\times \cdots \times n_d}$: the $d$-th order hypermatrices of dimensions $n_1,n_2,\cdots,n_d$.
		\item $[a,b]$: the set of integers $a\leq i \leq b$.
		\item $\langle b\rangle:=[1,b]$.
		\item $\ID(i_1,\cdots,i_d;n_1,\cdots,n_d)$: A set of ordered $d$ indices.
		%
		%\item $\lcm(n,p)$: the least common multiple of $n$ and $p$.
		%
		\item $\d_n^i$: the $i$-th column of the identity matrix $I_n$.
		\item $\d_n^I:=[(\d_n^1)^T, (\d_n^2)^T,\cdots,(\d_n^n)^T]^T$.
		\item $\D_n:=\left\{\d_n^i\vert i=1,\cdots,n\right\}$.
		
		\item $L\in \R^{m\times n}$ is a logical matrix, if its column set $\Col(L)\subset \D_m$. Denote by $\mathcal{L}_{m\times n}$ the set of $m\times n$ logical matrices.
		
		\item $\d_m[i_1,\cdots,i_n]:=[\d_m^{i_1},\cdots,\d_m^{i_n}]$.
		\item $\J_{\ell}:=(\underbrace{1,1,\cdots,1}_{\ell})^\mathrm{T}$.
		\item ${\bf V}_r(A)$: the row stacking form of $A$.
		\item ${\bf V}_c(A)$: the column stacking form of $A$.
		\item ${\bf V}^s_r(A)$: the s-row stacking form of $A$.
		\item ${\bf V}^s_c(A)$: the s-column stacking form of $A$.
		\item ${\bf S}_d$: the $d$-th order permutation group.
		%
		%\item $\cdet$: combinatorial hyperdeterminant.
		%
		%\item $\ddet$: modified combinatorial hyperdeterminant.
		%
		%\item $\Det$: slice-based hyperdeterminant.
		%
		\item $\ltimes$: the STP of matrices.
		\item $\times^{{\bf r}}_{{\bf s}}$: the contracted product of hypermatrices.
		\item $\ltimes^{{\bf r}}_{{\bf s}}$: the contracted STP of hypermatrices.
		%
		%\item $c^{{\bf r}}_{{\bf s}}$: contraction operator of hypermatrices.
		%
		%\item $\copyright^{{\bf r}}_{{\bf s}}$: contraction STP of hypermatrices.
		%
		%\item $\circledcirc$: STP of hypermatrices.
		%
		%\item $\GL(n^{(d)},\F)$: $n^{(d)}$-general linear group of hypercubics.
		%
		%\item $\GL(*^{(d)},\F)$: $d$-th order general linear group of hypercubics.
	\end{enumerate}

	\section{Matrix Expression of Hypermatrices}
	
	\subsection{Set-based Hypermatrix}
	
	\begin{dfn}\label{d2.1.1} \cite{lim13} For $n_1,\cdots,n_d\in \N$, a function $f:\langle n_1\rangle\times \cdots\times \langle n_d\rangle\ra \F$ is called a hypermatrix (over $\F$), or an order $d$ hypermatrix or a $d$-hypernatrix. A $d$-hypermatrix is also denoted by $A=[a_{i_1,\cdots,i_d}]$.
	\end{dfn}
	
	From a set point of view, a hypermatrix consists of two ingredients:
	\begin{itemize}
		\item[(i)] a set of data
		\begin{align}\label{2.1.1}
			D_A:=\{a_{i_1,\cdots,i_d}\;|\; i_s\in \langle n_s\rangle,\;s\in\langle d\rangle\};
		\end{align}
		\item[(ii)] an order over $D_A$, which makes $D_A$ a totally ordered set.
		In other words, each element in $D_A$ has a fixed position, just as in the case of a matrix.
	\end{itemize}
	
	To make the order of $D_A$ precise, we introduce a set of ordered indices, $\ID$, as follows:
	$$
	\ID(i_{\a_1},\cdots,i_{\a_d};n_{\a_1},\cdots,n_{a_d}),
	$$
	which stands for a set of data as in (\ref{2.1.1}), where $(\a_1,\cdots,\a_d)$ is a permutation of $\langle d\rangle$.  The order of $\ID$ is determined as follows:
	$a_{p_1,\cdots,p_d}\prec a_{q_1,\cdots,q_d}$ if and only if there is an $s\in \langle d\rangle$, such that $p_i=q_i$, $i<s$ and $p_s<q_s$.
	
	If $(\a_1,\cdots,\a_d)=(1,2,\cdots,d)$, the corresponding ordered set of indices is called a natural order, denoted by ${\bf i}_{\langle d\rangle}$.
	
	Based on the above argument, a hypermatrix can be defined alternatively as follows.
	
	\begin{dfn}\label{d2.1.2} For $n_1,\cdots,n_d\in \N$, a hypermatrix (over $\F$), is a set of order $d$ data as in (\ref{2.1.1}) with an $\ID$ as
		$\ID({\bf i}_{\langle d\rangle};n_1,\cdots,n_d)$.
	\end{dfn}
	
	To apply classical matrix analysis to hypermatrices, the matrix expression of hypermatrices is a key issue.
	
	Let ${\bf r}=(r_1,r_2,\cdots,r_d)$ be a set of indices, and
	\begin{align}\label{2.1.2}
		{\bf r}=(r_{i_1},r_{i_2},\cdots,r_{i_p})\cup (r_{j_1},r_{j_2},\cdots,r_{j_q}):={\bf r}_1\cup {\bf r}_2
	\end{align}
	be a partition, where ${\bf r}_1\cap {\bf r}_2=\emptyset$, and $p+q=d$.
	For ease of notation, we will assume that the elements in two subsets inherit the element order of the original set, unless otherwise stated.
	
	\begin{dfn}\label{d2.1.3} Given a hypermatrix $A=[a_{r_1,r_2,\cdots,r_d}]$. For each partition
		\begin{align}\label{2.1.3}
			{\bf r}_{\langle d\rangle}=\{r_{i_1},r_{i_2},\cdots,r_{i_p}\}\cup \{r_{j_1},r_{j_2},\cdots,r_{j_q}\}:={\bf r}_1\cup {\bf r}_2,
		\end{align}
		there is a matrix expression of $A$, denoted by
		\begin{align}\label{2.1.4}
			M_A^{{\bf r}_1\times {\bf r}_2}\in \F^{s\times t},
		\end{align}
		where
		$$
		s=\prod_{k=1}^pn_{i_k},\quad
		t=\prod_{\ell=1}^qn_{j_{\ell}}.
		$$
		Moreover, the elements in $M_A^{{\bf r}_1 \times {\bf r}_2}$ are $\{a_{r_1,r_2,\cdots,r_d}\}$, which are arranged
		by $\ID({\bf r}_1;n_{r_{i_1}},n_{r_{i_2}},\cdots,n_{r_{i_p}})$ for rows, and
		by $\ID({\bf r}_2;n_{r_{j_1}},n_{r_{j_2}},\cdots,n_{r_{j_q}})$ for columns.
	\end{dfn}
	
	\begin{rem}\label{r2.1.4}
		\begin{itemize}
			\item[(i)] $p$ is called the contra-variant order of $M_A^{{\bf r}_1\times {\bf r}_2}$ and $q$ is called the co-variant order of $M_A^{{\bf r}_1\times {\bf r}_2}$.
			\item[(ii)] If $p=0$, then we set $s=1$, and if $q=0$, then we set $t=1$.
			\item[(iii)] Particularly, in (\ref{2.1.3}) we require a natural order unless elsewhere is stated. That is,
			$$
			\begin{array}{l}
				r_{i_1}<r_{i_2}<\cdots<r_{i_p},\\
				r_{j_1}<r_{j_2}<\cdots<r_{j_q}.
			\end{array}
			$$
			Thus a $d$-hypermatrix has $2^d$ different matrix expressions with the natural order. The matrix expressions
			$M_A^{{\bf r}_1\times {\bf r}_2}$ can be briefly expressed as $M_A^{{\bf r}_1}$, if there is no possible confusion.

		\end{itemize}
	\end{rem}
	
	\begin{exa}\label{e2.1.5}
		Given $A=[a_{i_1 i_2 i_3}]\in \F^{2\times 3\times 2}$. Then
		\begin{itemize}
			\item[(i)]
			$$
			\begin{array}{ccl}
				M_A^{\emptyset}=[a_{111}~a_{112}~a_{121}~a_{122}~a_{131}~a_{132}~a_{211}~a_{212}~a_{221}~a_{222}~a_{231}~a_{232}];
			\end{array}
			$$
			\item[(ii)]
			$$
			M_A^{(1)}=
			\begin{bmatrix}
				a_{111}&a_{112}&a_{121}&a_{122}&a_{131}&a_{132}\\
				a_{211}&a_{212}&a_{221}&a_{222}&a_{231}&a_{232}
			\end{bmatrix};
			$$
			$$
			M_A^{(2)}=
			\begin{bmatrix}
				a_{111}&a_{112}&a_{211}&a_{212}\\
				a_{121}&a_{122}&a_{221}&a_{222}\\
				a_{131}&a_{132}&a_{231}&a_{232}
			\end{bmatrix}; \quad \mbox{etc.}
			$$
			\item[(iii)]
			$$
			M_A^{(1,2)}=
			\begin{bmatrix}
				a_{111}&a_{112}\\
				a_{121}&a_{122}\\
				a_{131}&a_{132}\\
				a_{211}&a_{212}\\
				a_{221}&a_{222}\\
				a_{231}&a_{232}
			\end{bmatrix};\quad \mbox{etc.}
			$$
			$$
			M_A^{(1,3)}=
			\begin{bmatrix}
				a_{111}&a_{121}&a_{131}\\
				a_{112}&a_{122}&a_{132}\\
				a_{211}&a_{221}&a_{231}\\
				a_{212}&a_{222}&a_{232}
			\end{bmatrix}; \quad \mbox{etc.}
			$$
			\item[(iv)]
			$$
			M_A^{(1,2,3)}=\left(M_A^{\emptyset}\right)^\mathrm{T}.
			$$
		\end{itemize}
	\end{exa}
	
	\begin{dfn}\label{d2.1.6}
		\begin{itemize}
			\item[(i)] $V_A:=M_A^{\emptyset}$
			is called the (row) vector expression of hypermatrix $A$.
			\item[(ii)]
			$M_A:=M_A^{(1)}$
			is called the contra-variant one matrix expression (briefly, matrix-1 expression) of hypermatrix $A$.
		\end{itemize}
	\end{dfn}
	
	\begin{dfn}\label{d2.3.1} Let $x_i\in \F^{n_i}$, $i\in\langle d\rangle$.
		\begin{itemize}
			\item[(i)]
			\begin{align}\label{2.3.1}
				x:=\ltimes_{i=1}^dx_i
			\end{align}
			is called a hypervector of order $d$.
			\item[(ii)] The set of order $d$ hypervectors with corresponding dimensions is denoted by $\F^{n_1\ltimes \cdots \ltimes n_d}$.
		\end{itemize}
	\end{dfn}
	
	Note that the components of $x$ can be expressed as
	$$
	\begin{array}{ccl}
		\com(x):=\left\{x_{i_1,\cdots,i_d}=x_1^{i_1}x_2^{i_2}\cdots x_d^{i_d}\;|\; i_j\in\langle n_j\rangle, j\in\langle d\rangle\right\},
	\end{array}
	$$
	where $x_r^{i_r}$ is the $i_r$-th component of $x_r$. Hence it is clear that $\com(x)$ (or briefly hypervector $x$)  is a hypermatrix of order $d$. Moreover, it is obvious that
	\begin{prp}\label{p2.3.2}
		\begin{align}\label{2.3.2}
			\F^{n_1\ltimes \cdots \ltimes n_d}\subset \F^{n_1\times \cdots \times n_d}
		\end{align}
		is a subset of hypermatrices.
	\end{prp}

	\subsection{$\sigma$-transpose of Hypermatrices}
	
	\begin{dfn}\label{d2.2.1} \cite{lim13}
		\begin{itemize}
			\item[(i)] Given a $d$-hypermatrix $A=\left[ a_{j_1,j_2,\cdots,j_d} \right]\in \F^{n_1\times n_2\times \cdots\times n_d}$, and assume $\sigma\in {\bf S}_d$.
			The $\sigma$-transpose of $A$ is
			\begin{align}\label{2.2.1}
				A^{\sigma}:=\left[ a_{j_{\sigma(1)}\cdots j_{\sigma(d)}}\right]\in \F^{n_{\sigma(1)}\times \cdots\times n_{\sigma(d)}}.
			\end{align}
			\item[(ii)] If a $d$-hypermatrix $A\in \F^{\overbrace{n\times \cdots\times n}^d}$, then $A$ is said to be $d$-hypercubic.
			\item[(iii)] A $d$-hypercubic $A\in \F^{\overbrace{n\times \cdots\times n}^d}$ is said to be symmetric if $\forall \sigma\in {\bf S}_d$, $A^{\sigma}=A$; $A$ is said to be skew-symmetric if $\forall \sigma\in {\bf S}_d$, $A^{\sigma}=\sgn(\sigma)A$.
		\end{itemize}
	\end{dfn}
	
	It is obvious that a matrix is a $2$-hypermatrix, so the above general definitions coincide with the corresponding definitions for matrices.
	
	\begin{prp}\label{p2.2.2} A
		$2$-hypercubic $A\in \F^{n\times n}$ is (skew-)symmetric, if and only if,
		$M_A$ is (skew-)symmetric.
	\end{prp}

	\begin{prp}\label{p2.2.3} Let $A\in \F^{n_1\times \cdots\times n_d}$ and ${\bf r}\subset {\bf d}=\langle d\rangle$. Then
		\begin{align}\label{2.2.100}
			\left[M_A^{{\bf r}\times ({\bf d}\backslash {\bf r})}\right]^T=M_A^{({\bf d}\backslash {\bf r})\times {\bf r}}.
		\end{align}
	\end{prp}
	
	\begin{rem}\label{r2.2.4} The above arguments stand true even when $\F$ is a set of perfect hypercomplex numbers (PHNs) \cite{che21b}. In fact, most of arguments throughout this paper also hold for PHNs.
	\end{rem}
	
	Next, we recall the permutation matrix \cite{che20}, which is a generalization of swap matrix.
	
	\begin{dfn}\label{d2.2.4} Let $n=\prod_{i=1}^dn_i$, $n_i\geq 2$, $\sigma\in {\bf S}_n$. A logical matrix $W_{[n_1,n_2,\cdots,n_d]}^{\sigma}\in {\cal L}_{n\times n}$, called a
		$\sigma$-permutation matrix, is constructed as follows:
		\begin{itemize}
			\item Step 1: Define
			$$
			\begin{array}{l}
				D=D_{[n_1,n_2,\cdots,n_d]}^{\sigma}:=\left\{\d^{j_1}_{n_{\sigma(1)}}\d^{j_2}_{n_{\sigma(2)}}\cdots \d^{j_d}_{n_{\sigma(d)}}\;\Big|\;j_i\in\langle n_{\sigma(i)}\rangle, i=1,2,\cdots,d
				\right\}.
			\end{array}
			$$
			\item Step 2: Arrange $\{\sigma(i)\;|\;i\in[1,d]\}$ into an increasing sequence as
			$$
			1=\sigma(i_1)<\sigma(i_2)<\cdots<\sigma(i_d)=d.
			$$
			That is,
			$$
			i_j=\sigma^{-1}(j),\quad j\in[1,d].
			$$
			Set an index order as
			$$
			\begin{array}{l}
				\ID:=\ID \left(j_{i_1},j_{i_2},\cdots,j_{i_d};n_{\sigma(i_1)},n_{\sigma(i_2)},\cdots,n_{\sigma(i_d)}\right)\\
				~=\ID \left(j_{\sigma^{-1}(1)},j_{\sigma^{-1}(2)},\cdots,j_{\sigma^{-1}(d)};n_{1},n_{2},\cdots,n_{d}\right).
			\end{array}
			$$
			\item Step 3:
			\begin{align}\label{2.2.2}%201.2}
			\begin{array}{l}
    W_{[n_1,n_2,\cdots,n_d]}^{\sigma}:=\left[\d^{j_1}_{n_{\sigma(1)}}\d^{j_2}_{n_{\sigma(2)}}\cdots \d^{j_d}_{n_{\sigma(d)}}\;\Big|~\mbox{arranged by the order of}~ \ID\right].
			\end{array}
		\end{align}
	\end{itemize}
\end{dfn}

\begin{exa}\label{e2.2.5} Consider $d=3$, $n_1=2$, $n_2=3$, and $n_3=5$. We construct $W^{\sigma}:=W^{\sigma}_{[2,3,5]}$.
\begin{itemize}
	\item[(1)] $\sigma_1=identity$ (i.e., $[1,2,3]\ra [1,2,3]$):
	We have $W^{\sigma_1}=I_{30}$.
	\item[(2)] $\sigma_2=(2,3)$ (i.e., $[1,2,3]\ra [1,3,2]$):
	Then
	$$
	D=\{\d_2^{j_1}\d_5^{j_2}\d_3^{j_3}\;|\;
	j_1\in[1,2],j_2\in[1,5],j_3\in[1,3]\}
	$$
	$$
	\begin{array}{l}
		W^{\sigma_2}=\left[\d_2^1\d_5^1\d_3^1,\d_2^1\d_5^2\d_3^1,\d_2^1\d_5^3\d_3^1,\d_2^1\d_5^4\d_3^1,\d_2^1\d_5^5\d_3^1,\d_2^1\d_5^1\d_3^2,\d_2^1\d_5^2\d_3^2,\d_2^1\d_5^3\d_3^2,\d_2^1\d_5^4\d_3^2,\d_2^1\d_5^5\d_3^2,\right.\\
		~~\d_2^1\d_5^1\d_3^3,\d_2^1\d_5^2\d_3^3,\d_2^1\d_5^3\d_3^3,\d_2^1\d_5^4\d_3^3,\d_2^1\d_5^5\d_3^3,\d_2^2\d_5^1\d_3^1,\d_2^2\d_5^2\d_3^1,\d_2^2\d_5^3\d_3^1,\d_2^2\d_5^4\d_3^1,\d_2^2\d_5^5\d_3^1,\\
		~~\left.\d_2^2\d_5^1\d_3^2,\d_2^2\d_5^2\d_3^2,\d_2^2\d_5^3\d_3^2,\d_2^2\d_5^4\d_3^2,\d_2^2\d_5^5\d_3^2,\d_2^2\d_5^1\d_3^3,\d_2^2\d_5^2\d_3^3,\d_2^2\d_5^3\d_3^3,\d_2^2\d_5^4\d_3^3,\d_2^2\d_5^5\d_3^3\right]\\
		~~=\d_{30}[ 1, 4, 7,10,13, 2, 5, 8,11,14, 3, 6, 9,12,15,\\
		16,19,22,25,28,17,20,23,26,29,18,21,24,27,30].
	\end{array}
	$$
	\item[(3)] $\sigma_3=(1,2)$ (i.e., $[1,2,3]\ra [2,1,3]$):
	
	Similarly, we have
	$$
	D=\{\d_3^{j_1}\d_2^{j_2}\d_5^{j_3}\;|\;
	j_1\in[1,3],j_2\in[1,2],j_3\in[1,5]\}
	$$
	$$
	\begin{array}{l}
		W^{\sigma_3}=\left[\d_3^1\d_2^1\d_5^1,\d_3^1\d_2^1\d_5^2,\cdots,\d_3^2\d_2^1\d_5^1,\cdots, \d_3^2\d_2^1\d_5^5,\cdots, \d_3^3\d_2^2\d_5^5\right]\\
		~~=\d_{30}[ 1, 2, 3, 4, 5,11,12,13,14,15,21,22,23,24,\\
		25,6, 7, 8, 9,10,16,17,18,19,20,26,27,28,29,30].
	\end{array}
	$$
	\item[(4)] $\sigma_4=(1,2,3)$ (i.e., $[1,2,3]\ra [2,3,1]$):
	
	We have
	$$
	\begin{array}{l}
		W^{\sigma_4}=\left[\d_3^1\d_5^1\d_2^1,\cdots,\d_3^3\d_5^5\d_2^2\right]\\
		~~=\d_{30}[ 1, 3, 5, 7, 9,11,13,15,17,19,21,23,25,27,\\
		29,2, 4, 6, 8,10,12,14,16,18,20,22,24,26,28,30].
	\end{array}
	$$
	\item[(5)] $\sigma_5=(1,3,2)$ (i.e., $[1,2,3]\ra [3,1,2]$):
	
	Then
	$$
	\begin{array}{l}
		W^{\sigma_5}=\left[\d_5^1\d_2^1\d_3^1,\cdots,\d_5^5\d_2^2\d_3^3\right]\\
		~=\d_{30}[ 1, 7,13,19,25, 2, 8,14,20,26, 3, 9,15,21,27\\
		4,10,16,22,28, 5,11,17,23,29, 6,12,18,24,30].
	\end{array}
	$$
	\item[(6)] $\sigma_6=(1,3)$ (i.e., $[1,2,3]\ra [3,2,1]$):
	
	Then
	$$
	\begin{array}{l}
		W^{\sigma_6}=\left[\d_5^1\d_3^1\d_2^1,\cdots,\d_5^5\d_3^3\d_2^2\right]\\
		~~=\d_{30}[ 1, 7,13,19,25, 3, 9,15,21,27, 5,11,17,23,\\
		29, 2, 8,14,20,26, 4,10,16,22,28, 6,12,18,24,30].
	\end{array}
	$$
	
\end{itemize}
\end{exa}

When $n_1=n_2=\cdots=n_d:=n$ the corresponding permutation is briefly denoted by
$$
W^{\sigma}_n:=W^{\sigma}_{[n_1,n_2,\cdots,n_d]}.
$$
Since this kind of permutation matrices are of particular importance, some of them are listed in the Appendix.

Some basic properties of permutation matrices are presented in the following proposition, which follows from definition immediately.

\begin{prp}\label{p2.2.6}
\begin{itemize}
	\item[(i)]
	\begin{align}\label{2.2.3}
		\left[W^{\sigma}_{[n_1,\cdots,n_d]}\right]^T=\left[W^{\sigma}_{[n_1,\cdots,n_d]}\right]^{-1}=
		W^{\sigma^{-1}}_{[n_1,\cdots,n_d]}.
	\end{align}
	\item[(ii)] Let $\sigma,\mu\in {\bf S}_d$. Then
	\begin{align}\label{2.2.4}
		W^{\sigma}_{[n_1,n_2,\cdots,n_d]}W^{\mu}_{[n_1,n_2,\cdots,n_d]}=
		W^{\sigma\circ \mu}_{[n_1,n_2,\cdots,n_d]}.
	\end{align}
\end{itemize}
\end{prp}

The following proposition shows the basic function of permutation matrices.

\begin{prp}\label{p2.2.7} \cite{che20} Assume $x_i\in \F^{n_i}$, $i\in \langle d\rangle$, $\sigma\in {\bf S}_d$. Then
\begin{align}\label{2.2.5}
	\ltimes_{i=1}^dx_{\sigma(i)}=W^{\sigma}_{[n_1,n_2,\cdots,n_d]}\ltimes_{i=1}^dx_{i}.
\end{align}
\end{prp}

As an immediate consequence we have the following result, which shows how to calculate $A^{\sigma}$.

\begin{prp}\label{p2.2.8} Let $A\in \F^{n_1\times \cdots \times n_d}$ be a hypermatrix of order $d$. Then
\begin{align}\label{2.2.6}
	V_{A^{\sigma}}=V_A \left[W_{[n_1,\cdots,n_d]}^{\sigma}\right]^T=V_A W_{[n_1,\cdots,n_d]}^{\sigma^{-1}}.
\end{align}
\end{prp}

\noindent{\it Proof.} Proposition \ref{p2.2.7} implies that
$$
V^T_{A^{\sigma}}=W_{[n_1,\cdots,n_d]}^{\sigma} V_A.
$$
Taking transpose on both sides yields (\ref{2.2.6}).
\hfill $\Box$

\subsection{Conversion of Matrix Expressions}

\begin{dfn}\label{d2.3.1}
\begin{itemize}
	\item[(i)]  Let $A=[a_{i,j}]\in \F^{m\times n}$ be a matrix. Then
	\begin{align}\label{2.3.1}
		{\bf V}_r(A):=[a_{1,1},a_{1,2},\cdots,a_{1,n},a_{2,1},\cdots,a_{m,n}]
	\end{align}
	is called the row stacking form of $A$;
	
	\begin{align}\label{2.3.2}
		{\bf V}_c(A):=[a_{1,1},a_{2,1},\cdots,a_{m,1},a_{1,2},\cdots,a_{m,n}]
	\end{align}
	is called the column stacking form of $A$.
	
	\item[(ii)] Let $x\in \F^n$ and $s|n$. Say, $n=st$. Then
	\begin{align}\label{2.3.3}
		{\bf V}^s_r(x):=
		\begin{bmatrix}
			x_1&x_2&\cdots&x_s\\
			x_{s+1}&x_{s+2}&\cdots&x_{2s}\\
			\vdots&~&\ddots&\vdots\\
			x_{(t-1)s+1}&x_{(t-1)s+2}&\cdots&x_{ts}\\
		\end{bmatrix}
	\end{align}
	\begin{align}\label{2.3.4}
		{\bf V}^s_c(x):=
		\begin{bmatrix}
			x_1&x_{s+1}&\cdots&x_{(t-1)s+1}\\
			x_2&x_{s+2}&\cdots&x_{(t-1)s+2}\\
			\vdots&~&\ddots&\vdots\\
			x_{s}&x_{2s}&\cdots&x_{ts}\\
		\end{bmatrix}
	\end{align}
	
	\item[(iii)] Let $A\in \F^{m\times n}$ and $s|(mn)$. Then
	\begin{align}\label{2.3.5}
		{\bf V}^s_r(A):={\bf V}^s_r\left({\bf V}_r(A)\right)
	\end{align}
	is called the $s$-row stacking form.
	\begin{align}\label{2.3.6}
		{\bf V}^s_c(A):={\bf V}^s_c\left({\bf V}_c(A)\right)
	\end{align}
	is called the $s$-column stacking form.
\end{itemize}
\end{dfn}

\begin{prp}\label{p2.3.2} \cite{che07} Let $A\in \F^{m\times n}$, $X\in \F^{n\times q}$, and $Y\in \F^{p\times m}$. Then
\begin{align}\label{2.3.7}
	{\bf V}_r(AX)=A\ltimes {\bf V}_r(X),
\end{align}
\begin{align}\label{2.3.8}
	{\bf V}_c(YA)=A^T\ltimes {\bf V}_c(Y).
\end{align}
\end{prp}

Denote by
$$
\d_n^I:={\bf V}_r(I_n)={\bf V}_c(I_n)=[(\d_n^1)^T, (\d_n^2)^T,\cdots,(\d_n^n)^T]^T.
$$

\begin{prp}\label{p2.3.3} Let $A\in \F^{m\times n}$. Then
\begin{align}\label{2.3.9}
	{\bf V}_r(A)=A\ltimes \d_n^I.
\end{align}

\begin{align}\label{2.3.10}
	{\bf V}_c(A)=A^T\ltimes \d_m^I.
\end{align}

Conversely,

\begin{align}\label{2.3.11}
	A={\bf V}_r^n({\bf V}_r(A))={\bf V}_c^m({\bf V}_c(A)).
\end{align}
\end{prp}

\noindent{\it Proof.} (\ref{2.3.9}) and (\ref{2.3.10}) come from Proposition \ref{p2.3.2} immediately.
(\ref{2.3.11}) follows from the definition.
\hfill $\Box$

By definition and (\ref{2.3.11}), we have
\begin{align}\label{2.3.1101}
A={\bf V}_r^n(A)={\bf V}_c^m(A).
\end{align}

Set ${\bf i_r}=(i_1,\cdots,i_r)\subset {\bf d}=\langle d\rangle$, and denote
\begin{align}\label{2.3.12}
\sigma_{{\bf i_r}}: {\bf d}\ra ({\bf i_r},{\bf d}\backslash {\bf i_r}).
\end{align}

The following proposition shows how to convert a matrix expression to a vector form and vise versa for  a hypermatrix.

\begin{prp}\label{p2.3.4} Given $A=[a_{i_1,\cdots,i_d}]\in \F^{n_1\times\cdots\times n_d}$,
${\bf i_r}=(i_1,\cdots,i_r)\subset {\bf d}=\langle d\rangle$, and $\sigma_{{\bf i_r}}$ is as in (\ref{2.3.12}),
$$
n_{{\bf i_r}}=\prod_{s=1}^rn_{i_s},\quad n_{{\bf d}\backslash {\bf i_r}}=\prod_{i_j\in {\bf d}\backslash {\bf i_r}}n_{i_j}.
$$
Then
\begin{itemize}
	\item[(i)] (Vector Form to Matrix Form)
	\begin{align}\label{2.3.13}
		M_A^{{\bf i_r}\times ({\bf d}\backslash {\bf i_r})}=
		{\bf V}_r^{n_{{\bf d}\backslash {\bf i_r}}}\left(V_A W^{\sigma^{-1}_{{\bf i_r}}}_{[n_1,\cdots,n_d]}\right).
	\end{align}
	
	\item[(ii)]  (Matrix Form to Vector Form)
	\begin{align}\label{2.3.14}
		V_A=\left(M_A^{{\bf i_r}\times ({\bf d}\backslash {\bf i_r})}\ltimes \d^I_{n_{{\bf d}\backslash {\bf i_r}}}\right)^T
		W_{[n_1,\cdots,n_d]}^{\sigma_{{\bf i_r}}}.
	\end{align}
\end{itemize}
\end{prp}

\noindent{\it Proof.} For (i), we have
$$
\begin{array}{l}
	M_A^{{\bf i_r}\times ({\bf d}\backslash {\bf i_r})}=
	{\bf V}_r^{n_{{\bf d}\backslash {\bf i_r}}}(V_{A^{\sigma_{\bf i_r}}})={\bf V}_r^{n_{{\bf d}\backslash {\bf i_r}}}\left(V_A W^{\sigma^{-1}_{{\bf i_r}}}_{[n_1,\cdots,n_d]}\right).
\end{array}
$$
As for (ii), using (\ref{2.2.6}) and (\ref{2.3.9}), we have
$$
\begin{array}{l}
	V_A=V_{A^{\sigma_{{\bf i_r}}}}W_{[n_1,\cdots,n_d]}^{\sigma_{{\bf i_r}}}=\left(M_A^{{\bf i_r}\times ({\bf d}\backslash {\bf i_r})}\ltimes \d^I_{n_{{\bf d}\backslash {\bf i_r}}}\right)^T
	W_{[n_1,\cdots,n_d]}^{\sigma_{{\bf i_r}}}.
\end{array}
$$
\hfill $\Box$

Using (\ref{2.3.13}) and (\ref{2.3.14}), the formula transforming one matrix form to another can be obtained.

\begin{cor}\label{c2.3.5} Let ${\bf i_r}$, $\sigma_{{\bf i_r}}$ be as in Proposition \ref{p2.3.4}, and
${\bf j_s}=(j_1,\cdots,j_s)$ and $\sigma_{{\bf j_s}}: {\bf d}\ra ({\bf j_s}, {\bf d}\backslash {\bf j_s})$. Then
\begin{align}\label{2.3.15}%{2.3.9}
	\begin{array}{ccl}
		M_A^{{\bf j_s}\times ({\bf d}\backslash {\bf j_s})}
		&=&{\bf V}_r^{n_{{\bf d}\backslash {\bf j_s}}}\left[
		\left(
		M_A^{{\bf i_r}\times ({\bf d}\backslash {\bf i_r})}\ltimes \d_{{\bf d}\backslash {\bf i_r}}^I\right)^T\times W^{\sigma_{{\bf i_r}}}_{[n_1,\cdots,n_d]} W^{\sigma^{-1}_{{\bf j_s}}}_{[n_1,\cdots,n_d]}
		\right].
	\end{array}
\end{align}
\end{cor}

\section{Contracted Product of Hypermatrices}

\begin{dfn}\label{d3.1.1} \cite{lim13} Let $A=[a_{i_1,\cdots,i_d}]\in \F^{n_1\times \cdots\times n_{d-1}\times n}$, $B=[b_{j_1,\cdots,j_r}]\in \F^{n\times m_2\times \cdots \times m_{r}}$. Then the contracted product of $A$ and $B$, denoted by
	\begin{align}\label{3.1.1}
		C=A \times^{i_d}_{j_1} B,
	\end{align}
	is s hypermatrix $C=[c_{i_1,\cdots,i_{d-1},j_2,\cdots,j_r}]\in \F^{n_1\times \cdots\times n_{d-1}\times m_2\times \cdots \times m_r}$, where
	\begin{align}\label{3.1.2}
		c_{i_1,\cdots,i_{d-1},j_2,\cdots,j_r}=\dsum_{k=1}^n a_{i_1,\cdots,i_{d-1},k}b_{k,j_2,\cdots,j_r}.
	\end{align}
\end{dfn}

\begin{prp}\label{p3.1.2}Let $M_A^{(i_1,\cdots,i_{d-1})\times i_d}$ and $M_B^{j_1\times (j_2,\cdots,j_r)}$ be the corresponding matrix expressions of $A$ and $B$ respectively. Then
	\begin{align}\label{3.1.3}
		M_C^{(i_1,\cdots,i_{d-1})\times(j_2,\cdots,j_r)}=M_A^{(i_1,\cdots,i_{d-1})\times i_d}\times M_B^{j_1\times (j_2,\cdots,j_r)},
	\end{align}
	where the $\times$ between $M_A$ and $M_B$ is the conventional matrix product.
\end{prp}

Definition \ref{d3.1.1} can be extended to the case of multiple common indices.

\begin{dfn}\label{d3.1.3} Let $A=[a_{i_1,\cdots,i_d}]\in \F^{n_1\times \cdots\times n_d}$, $B=[b_{j_1,\cdots,j_r}]\in \F^{m_1\times \cdots\times m_{r}}$ with
	$$
	n_{\a_t}=m_{\b_t}:=\ell_t,\quad \a_t\in \langle d\rangle, \b_t\in \langle r\rangle, t\in \langle s\rangle,
	$$
	where $s\leq \min(d,r)$ is the number of equal dimension indexes. Define
	\begin{align}\label{3.1.4}
		\begin{array}{l}
			A\times^{i_{\a_1},\cdots,i_{\a_s}}_{j_{\b_1},\cdots,j_{\b_s}} B:=C\\
			~~\in \F^{n_1\times \cdots \times \hat{n}_{\a_1}\times \cdots\times \hat{n}_{\a_s}\times \cdots\times n_d\times
				m_1\times \cdots\times \hat{m}_{\b_1}\times \cdots\times \hat{m}_{\b_s}\times \cdots\times m_r},
		\end{array}
	\end{align}
	where
	$$
	\begin{array}{l}
		c_{i_1,\cdots,\hat{i}_{\a_1},\cdots,\hat{i}_{\a_s},\cdots,i_d,
			j_1,\cdots,\hat{j}_{\b_1},\cdots,\hat{j}_{\b_s},\cdots,j_r}\\
		=\dsum_{k_1=1}^{\ell_1}\cdots \dsum_{k_s=1}^{\ell_s}a_{i_1\cdots k_1\cdots k_s\cdots i_d}b_{j_1\cdots k_1\cdots k_s\cdots j_d},
	\end{array}
	$$
	where a caret over any entry means that the respective entry is omitted.
\end{dfn}

\begin{prp}\label{p3.1.4}
	\begin{align}\label{3.1.5}
		\begin{array}{l}
			M_C^{(i_1,\cdots,\hat{i}_{\a_1},\cdots,\hat{i}_{\a_s},\cdots,i_d)\times
				(j_1,\cdots,\hat{j}_{\b_1},\cdots,\hat{j}_{\b_s},\cdots,j_r)}\\
			=M_A^{(i_1,\cdots,\hat{i}_{\a_1},\cdots,\hat{i}_{\a_s},\cdots,i_d)\times
				(i_{\a_1},\cdots,i_{\a_s})}\\
			~~\times M_B^{(j_{\b_1},\cdots,j_{\b_s})\times (j_1,\cdots,\hat{j}_{\b_1},\cdots,\hat{j}_{\b_s},\cdots,j_r)
			},
		\end{array}
	\end{align}
	where the $\times$ between $M_A$ and $M_B$ is the conventional matrix product.
\end{prp}

\begin{exa}\label{e3.1.5}
	Assume
	$A=[a_{i_1,i_2,i_3}]\in \F^{2\times 3\times 4}$, $B=[b_{j_1,j_2,j_3}]\in \F^{4\times 5\times 3}$  with natural ID, calculate
	$$
	C=A\times^{(i_2,i_3)}_{(j_3,j_1)} B.
	$$
	
	We have
	$$
	\begin{array}{l}
		M_A^{i_1\times (i_2,i_3)}=
		\left[
		\begin{array}{llllllllllll}
	a_{111}&a_{112}&a_{113}&a_{114}&a_{121}&a_{122}&a_{211}&a_{212}&a_{213}&a_{214}&a_{221}&a_{222}\\
	a_{123}&a_{124}&a_{131}&a_{132}&a_{133}&a_{134}&a_{223}&a_{224}&a_{231}&a_{232}&a_{233}&a_{234}\\
		\end{array}\right]\\
	\end{array}
	$$
	$$
	\begin{array}{l}
		M_B^{(j_3,j_1)\times j_2}=
		\begin{bmatrix}
			b_{111}&b_{121}&b_{131}&b_{141}&b_{151}\\
			b_{211}&b_{221}&b_{231}&b_{241}&b_{251}\\
			b_{311}&b_{321}&b_{331}&b_{341}&b_{351}\\
			b_{411}&b_{421}&b_{431}&b_{441}&b_{451}\\
			b_{112}&b_{122}&b_{132}&b_{142}&b_{152}\\
			b_{212}&b_{222}&b_{232}&b_{242}&b_{252}\\
			b_{312}&b_{322}&b_{332}&b_{342}&b_{352}\\
			b_{412}&b_{422}&b_{432}&b_{442}&b_{452}\\
			b_{113}&b_{123}&b_{133}&b_{143}&b_{153}\\
			b_{213}&b_{223}&b_{233}&b_{243}&b_{253}\\
			b_{313}&b_{323}&b_{333}&b_{343}&b_{353}\\
			b_{413}&b_{423}&b_{433}&b_{443}&b_{453}\\
		\end{bmatrix}
	\end{array}
	$$
	Then $C\in \F^{2\times 5}$ with
	$$
	M_C^{i_1\times j_2}=M_A^{i_1\times (i_2,i_3)}M_B^{(j_3,j_1)\times j_2}.
	$$
\end{exa}

A particular case is $A=[a_{i_1,\cdots,i_d}]\in \F^{n_1\times \cdots\times n_d}$, $B=[b_{i_{r_1},\cdots,i_{r_s}}]\in \F^{n_{r_1}\times \cdots\times n_{r_s}}$, and ${\bf r_s}:=\{r_1,\cdots,r_s\}\subset {\bf d}:=\langle d\rangle$. Then we briefly denote
$$
A\times_{{\bf r_s}}B:=A\times^{{\bf r_s}}_{\bf r_s}B.
$$
We call this kind of product the onto contracted product. They are of particular importance.

In the onto contracted product, $A$ can be considered as a multilinear mapping from
$\F^{n_{r_1}\times \cdots \times n_{r_s}}$ to $\F^{n_1\times \cdots \times \hat{n}_{r_1}\times \cdots \times \hat{n}_{r_s}\times \cdots \times n_d}$.
\begin{prp}\label{p3.1.6}
	Assume $A\in \F^{n_{1}\times \cdots \times n_{d}}$ and $B\in \F^{n_{r_1}\times \cdots \times n_{r_s}}$, where ${\bf r_s}:=\{r_1,\cdots,r_s\}\subset {\bf d}$. The contracted product of $A$ and $B$, denoted by
	$$
	C=A\times_{{\bf r_s}} B,
	$$
	can be obtained by one of the following two equivalent formulas:
	\begin{itemize}
		\item[(i)] Let
		$$
		M^A:=M_A^{({\bf d}\backslash {\bf r_s})\times {\bf r_s}}.
		$$
		Then
		\begin{align}\label{3.1.6}
			V_C^T=M^AV_B^T.
		\end{align}
		\item[(ii)] Let
		$\sigma\in {\bf S}_d$ be the permutation
		$$
		\sigma: {\bf d} \ra ( {\bf d}\backslash {\bf r_s}, \bf{r_s}).
		$$
		Then
		\begin{align}\label{3.1.7}
			V_C=V_A W^{\sigma^{-1}}V_B^T.
		\end{align}
	\end{itemize}
\end{prp}

\begin{dfn}\label{d3.1.7}
	\begin{itemize}
		\item[(i)] Let $A\in \F^{n_1\times \cdots\times n_{d}\times n_{d+1}\times\cdots\times n_{2d}}$, $B\in  \F^{n_1\times \cdots\times n_{d}}$, where
		$n_{d+i}=n_i$, $i\in\langle d\rangle$. Then
		$
		A:  \F^{n_1\times \cdots\times n_{d}}\ra \F^{n_1\times \cdots\times n_{d}}
		$
		is defined by
		\begin{align}\label{3.1.8}
			A(B):=A\times^{d+1,\cdots,2d}_{1,2,\cdots,d}B \in \F^{n_1\times \cdots\times n_{d}}.
		\end{align}
		This contracted product is called a unary operator on $\F^{n_1\times \cdots\times n_{d}}$.
		\item[(ii)] Let $A\in \F^{n_1\times \cdots\times n_{3d}}$, $B,C\in  \F^{n_1\times \cdots\times n_{d}}$, where
		$n_{2d+i}=n_{d+i}=n_i$, $i\in\langle d\rangle$. Then
		$
		A:  \F^{n_1\times \cdots\times n_{d}}\times \F^{n_1\times \cdots\times n_{d}} \ra \F^{n_1\times \cdots\times n_{d}}
		$
		is defined by
		\begin{align}\label{3.1.9}
			\begin{array}{ccl}
				A(B,C):=[A\times^{2d+1,\cdots,3d}_{1,2,\cdots,d}B]\times^{d+1,\cdots,2d}_{1,\cdots,d} C\in\F^{n_1\times \cdots\times n_{d}}.
		\end{array}\end{align}
		This contracted product is called a binary operator on $\F^{n_1\times \cdots\times n_{d}}$.
		\item[(iii)] Similarly, we can define $k$ argument (i.e., $k$ hypermatrix) operators for $k\geq 3$.
	\end{itemize}
\end{dfn}

\section{STP Realization of Contracted Product of Hypermatrices}

\begin{figure}
	\centering
	\setlength{\unitlength}{5 mm}
	\begin{picture}(17,11)(-0.5,-0.5)
		%\put(1,1){\framebox(4,13){}}
		%\put(6,1){\framebox(8,13){}}
		\thicklines
		\put(0,8){\framebox(6,2){Hypermatrix}}
		\put(-0.25,4){\framebox(6.5,2){Contracted Product}}
		\put(-0.75,0){\framebox(7.5,2){Resulting Hypermatrix}}
		\put(10,8){\framebox(6,2){Matrix}}
		\put(10,4){\framebox(6,2){STP}}
		\put(10,0){\framebox(6,2){Resulting Matrix}}
		\put(6,9){\vector(1,0){4}}
		\put(10,1){\vector(-1,0){3.25}}
		\put(3,8){\vector(0,-1){2}}
		\put(13,8){\vector(0,-1){2}}
		\put(3,4){\vector(0,-1){2}}
		\put(13,4){\vector(0,-1){2}}
		\put(7.1,9.2){Matrix}
		\put(6.5,8.2){Expression}
		\put(7,1.2){Converse}
		%\put(7.2,0.2){Back}
	\end{picture}
	\caption{STP-Realization of Contracted Product \label{Fig.4.1}}
\end{figure}
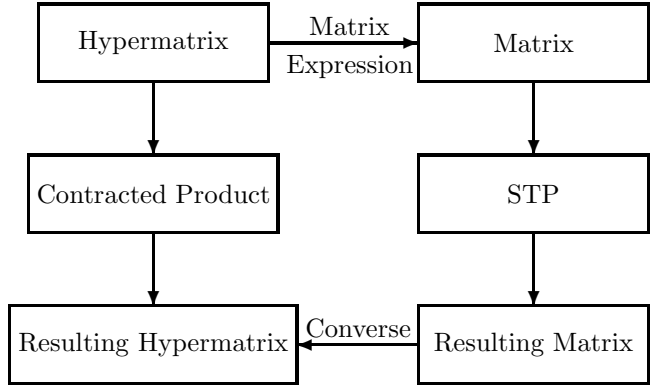

This section shows that a fundamental role of STP is to realize the contracted product of hypermatrices. This factor may help us to understand what is the essential meaning of STP. For this purpose, the hypermatrices are first expressed in their matrix expressions. The STPs are operators on matrices. When the hypermatrices are transformed into their matrix forms, the action of operators (especially certain contracted products) on their objective hypermatrices can be realized by the action of STPs on the matrix expressions of the corresponding hypermatrices. Figure \ref{Fig.4.1} shows this process.

Since there are many possible multilinear operators, including contracted products, over different hypermatrices, the corresponding operators over their matrix expressions are also diverse. To express these operators, the STPs are also diverse.

\subsection{Classical STPs}

From linear algebra one sees that the matrix product has two fundamental types: (1) Matrix-Matrix Product: It represents the composition of two linear mappings. (2) Matrix-Vector Product: It realizes a linear mapping over a vector space (or between two vector spaces). Fortunately,
when the dimension-matching condition is satisfied, the classical matrix product can realize these two functions simultaneously.

When the conventional matrix product is extended to STP, where the dimension-matching condition is not satisfied, an STP cannot realize these two functions simultaneously. Therefore, we need to distinguish between M-M STP and M-V products. First, we review the classical STPs \cite{che19b}:

\begin{dfn}\label{d3.2.1}
	\begin{itemize}
		\item[(i)] Let $A\in \F^{m\times n}$, $B\in \F^{p\times q}$ and $t=\lcm(n,p)$. The matrix-matrix (M-M) STP of $A$ and $B$ is defined as
		\begin{align}\label{3.2.1}
			A\ltimes B:=\left(A\otimes I_{t/n}\right)\left(B\otimes I_{t/p}\right).
		\end{align}
		\item[(ii)] Let $A\in \F^{m\times n}$, $x\in \F^{p}$ and $t=\lcm(n,p)$. The matrix-vector (M-V) STP of $A$ and $x$ is defined as
		\begin{align}\label{3.2.101}
			A\lvtimes x:=\left(A\otimes I_{t/n}\right)\left(x\otimes \J_{t/p}\right).
		\end{align}
		
		\item[(iii)] Let $x\in \F^{m}$, $y\in \F^{n}$ and $t=\lcm(m,n)$. The vector-vector (V-V) STP of $x$ and $y$ is defined as
		\begin{align}\label{3.2.102}
			x\vec{\cdot}y:=(x\otimes \J_{t/m})^T(y\otimes \J_{t/n}) \in \F,
		\end{align}
	\end{itemize}
\end{dfn}

Define
$$
{\cal M}=\bigcup_{m=1}^{\infty}\bigcup_{n=1}^{\infty}{\cal M}_{m\times n},\quad\R^{\infty}=\bigcup_{n=1}^{\infty}\R^n.
$$

Then the M-M STP can be considered as the product over ${\cal M}$; the M-V STP can be considered as the action of ${\cal M}$ on $\R^{\infty}$; and the V-V STP can be considered as an inner product over $\R^{\infty}$. It is obvious that they are the generalizations of the corresponding matrices with matrix-matrix or matrix-vector product, or vectors with vector product. Furthermore, they satisfy the following general properties.

\begin{prp}\label{p3.2.101}
	\begin{itemize}
		\item[(i)] (Associativity)
		$$
		(A \ltimes B)\ltimes C=A\ltimes (B\ltimes C),\quad A,B,C\in {\cal M}.
		$$
		$$
		(A\ltimes B)\lvtimes x=A\lvtimes (B\lvtimes x),\quad x\in \R^{\infty}.
		$$
		%$$
		%(x~\vec{\cdot}~ y)~\vec{\cdot}~ z=x ~\vec{\cdot}~ (y~\vec{\cdot}~ z).
		%$$

		\item[(ii)] (Distributivity) For $A,B,C\in {\cal M}$, $x,y,z\in \R^{\infty},$
		$$
		\begin{array}{l}
			(A + B)\ltimes C=A\ltimes C+B\ltimes C,\\
			C\ltimes (A+B)=C\ltimes A+C\ltimes B;\\
			(A + B)\lvtimes x=A\lvtimes x+B\lvtimes x,\\
			A\lvtimes (x+y)=A\lvtimes x+A\lvtimes y;\\
			(x+y)~\vec{\cdot}~ z=x~\vec{\cdot}~ z+y~\vec{\cdot}~ z,\\
			z~\vec{\cdot}~ (x+y)=z~\vec{\cdot}~ x+z~\vec{\cdot}~ y.\\
		\end{array}
		$$
	\end{itemize}
\end{prp}

\begin{rem}\label{r3.2.102} M-V and V-V STPs are not so commonly used as M-M STP. We explain it a little bit more.
	
	\begin{itemize}
		\item Topology on $\R^{\infty}$:
	\end{itemize}
	
	Consider $x\in \R^p$, $y\in \R^q$, ($x,y\in \R^{\infty}$), $t=\lcm(p,q)$. Define
	$$
	x~\vec{\pm}~y:=\left(x\otimes \J_{t/p}\right)\pm\left(y\otimes \J_{t/q}\right),
	$$
	then $\R^{\infty}$ becomes a (pseudo-) vector space.
	
	Furthermore,
	we define
	\begin{itemize}
		\item[(i)] (Inner Product):
		$$
		\langle x,y\rangle_{{\cal V}}:=\frac{1}{t} (x~\vec{\cdot}~y).
		$$
		\item[(ii)] (Norm):
		$$
		\|x\|_{{\cal V}}:=\sqrt{\langle x,x\rangle_{{\cal V}}}.
		$$
		\item[(iii)] (Distance):
		$$
		d(x,y)=\|x~\vec{-}~y\|_{{\cal V}}.
		$$
	\end{itemize}
	With this distance $\R^{\infty}$ becomes a topological space \cite{che19b}.
	
	\begin{itemize}
		\item Linear Dynamic Systems over $\R^{\infty}$:
	\end{itemize}

	A linear dynamic system over $\R^{\infty}$ is defined by
	$$
	x(t+1)=A \vec{\ltimes} x(t).
	$$
	
	It is a cross-dimensional system \cite{chepr}.
\end{rem}

\subsection{STP for Vector Expression of Hypermatrices}

Consider a multilinear mapping $\pi:\prod_{i=1}^dV_i\ra \F$, where $V_i=\F^{n_i}$, $i\in\langle d\rangle$. Assume
$$
\begin{array}{l}
	\pi\left(\d_{n_1}^{i_1},\cdots,\d_{n_d}^{i_d}\right)=c_{i_1,\cdots,i_d},\quad i_s\in \langle n_s\rangle,\; s\in \langle d\rangle.
\end{array}
$$
Let $\Pi:=[c_{i_1,\cdots,i_d}]$
with $\ID=\ID(i_1,\cdots,i_d;n_1,\cdots,n_d)$ be a hypermatrix of order $d$. Then the following result is well known.

\begin{prp}\label{p3.2.2} The multilinear mapping $\pi$ can be calculated by
	\begin{align}\label{3.2.2}
		\pi(x_1,\cdots,x_d)=V_{\Pi}\ltimes_{i=1}^dx_i,\quad x_i\in V_i,~i\in \langle d\rangle.
	\end{align}
\end{prp}

\begin{exa}\label{e3.2.201} Consider a finite game of $n$ players \cite{che15}. Assume player $i$ has $k_i$ strategies, denoted by
	$\d_{k_i}^{j_i}$, $j_i\in \langle k_i\rangle$, $i\in \langle n\rangle$. Let
	$$
	c_i(\d_{k_1}^{j_1},\cdots,\d_{k_n}^{j_n})=r^i_{j_i,\cdots,j_n},\quad i\in\langle n\rangle.
	$$
	Then
	$$
	D_i=\{r^i_{j_i,\cdots,j_n}\}\in \R^{k_1\times \cdots k_n}.
	$$
	Finally, we have
	$$
	c_i=V_{D_i}\ltimes_{t=1}^nx_t,\quad i\in\langle n\rangle.
	$$
\end{exa}

\subsection{STP for Matrix-1 Expression of Hypermatrices}

Consider a multilinear mapping $\pi:\prod_{i=1}^dV_i\ra V_0$, where $V_i=\F^{n_i}$, $i\in [0,d]$. Assume
$$
\pi\left(\d_{n_1}^{i_1},\cdots,\d_{n_d}^{i_d}\right)=\dsum_{i_0=1}^{n_0}c^{i_0}_{i_1,\cdots,i_d}\d_{n_0}^{i_0},\quad i_s\in \langle n_s\rangle,\; s\in \langle d\rangle.
$$
Then $\Pi:=[c^{i_0}_{i_1,\cdots,i_d}]$
is a hypermatrix of order $d+1$ with matrix expression $M_{\Pi}^{\ID}$, where $\ID=\ID(i_0;n_0)\times \ID(i_1,\cdots,i_d;n_1,\cdots,n_d)$. Then the following result is well known.

\begin{prp}\label{p3.2.4} The multilinear mapping $\pi$ can be calculated by
	\begin{align}\label{3.2.5}
		\pi(x_1,\cdots,x_d)=M_{\Pi}^{\ID}\ltimes_{i=1}^dx_i,\quad x_i\in V_i,~i\in \langle d\rangle.
	\end{align}
\end{prp}

\begin{exa}\label{3.2.5}
	\begin{itemize}
		\item[(i)] Cross Product on $\R^3$:
		
		Consider the cross product on $R^3$, denoted by $\vec{\times}$. Denote by
		$$
		\d_3^i\vec{\times} \d_3^j=\dsum_{k=1}^3c_{i,j}^k,\quad i,j\in\langle 3\rangle.
		$$
		Then we have
		$$
		M_C^{3\times (1,2)}=
		\begin{bmatrix}
			0&0&0&0&0&1&0&-1&0\\
			0&0&-1&0&0&0&1&0&0\\
			0&1&0&-1&0&0&0&0&0\\
		\end{bmatrix}
		$$
		It follows that
		$$
		%\begin{array}{ccl}
		Z=X\vec{\times} Y=M_C^{3\times (1,2)}\ltimes X\ltimes Y,\quad
		X,Y\in \R^3.
		%\end{array}
		$$
		\item[(ii)] General Linear Algebra $\gl(2,\R)$.
		Denote by
		$$
		\begin{array}{l}
			A_1=\begin{bmatrix}
				1&0\\
				0&0\\
			\end{bmatrix}\sim \vec{A}_1=\d_4^1\\
			A_2=\begin{bmatrix}
				0&1\\
				0&0\\
			\end{bmatrix}\sim \vec{A}_2=\d_4^2\\
			A_3=\begin{bmatrix}
				0&0\\
				1&0\\
			\end{bmatrix}\sim \vec{A}_3=\d_4^3\\
			A_4=\begin{bmatrix}
				0&0\\
				0&1\\
			\end{bmatrix}\sim \vec{A}_4=\d_4^4\\
		\end{array}
		$$
		Then $\d_4^i$, $i\in\langle 4\rangle$ form a basis of $\gl(2,\R)$.
		Denote by
		$$
		[A_i,A_j]=\dsum_{k=1}^4c_{i,j}^kA_k,\quad i,j\in\langle 4\rangle.
		$$
		Then $C=[c_{i,j}^k]\in \R^{4\times 4\times 4}$.
		Express $C$ into matrix form
		$$
		\begin{array}{ccl}
			M_C^{k,(i.j)}&=&\d_{4}[0,2,-3,0,-2,0,1-4,0,\\
			~&~&3,-1+4,0,-3,0,0,3,0],
		\end{array}
		$$
		where
		$$
		\begin{array}{l}
			\d_4^0=[0,0,0,0]^T,\\
			\d_4^{-3}=[0,0,-1,0]^T,\\
			\d_4^{1-4}=[1,0,0,-4]^T,\quad etc.\\
		\end{array}
		$$
		Let $X\in \gl(2,R)$ with
		$\vec{X}={\bf V}_c(X)$.
		Then
		$$
		\overrightarrow{[X,Y]}=M_C^{k,(i.j)}\ltimes \vec{X}\ltimes \vec{Y}.
		$$
	\end{itemize}

\end{exa}

\begin{rem}\label{r3.2.6}
	\begin{itemize}
		\item[(i)] Roughly speaking, in classical sense, an STP of matrices is a multilinear operator over hypermatrices.
		To be precise, when the hypermatrices are expressed into their matrix forms, the STP works as a matrix product.
		\item[(ii)]
		Since multilinear mappings over hypermatrices can be various, there can also be various STPs.
		\item[(iii)] When partial arguments are known, an operator becomes a restricted operator over the remaining arguments. In this case, STP combined with permutation matrices becomes more powerful.
	\end{itemize}
\end{rem}

\subsection{STP for General Matrix Expression of Hypermatrices}

We give two examples for this.

\begin{exa}\label{e3.5.1}
	
	Assume $V$ is an $n$-dimensional vector space over $\F$. $T: V^r\times (V^*)^s\ra \F$ is a tensor of covariant order $r$ and contra-variant order $s$. Let $\d_n^i$, $i\in \langle n\rangle$ be a basis of $V$ and $\omega_i=(\d_n^i)^T$, $i\in \langle n\rangle$ be the dual basis of dual space $V^*$, and set
	$$
	\begin{array}{l}
		T(\d_n^{i_1},\cdots,\d_n^{i_r},(\d_n^{j_1})^*,\cdots,(\d_n^{j_s})^*)
		=\mu^{i_1,\cdots,i_r}_{j_1,\cdots,j_s},\\
		~~ i_p\in\langle n\rangle, p\in\langle r\rangle,j_q\in\langle n\rangle, q\in\langle s\rangle.
	\end{array}
	$$

	Then
	$$
	\begin{array}{ccl}
		\Omega&:=&\{\mu^{i_1,\cdots,i_r}_{j_1,\cdots,j_s}\;|\; i_p\in\langle n\rangle, p\in\langle r\rangle, j_q\in\langle n\rangle, q\in\langle s\rangle\}\in \R^{\overbrace{n\times \cdots\times n}^{r+s}}
	\end{array}
	$$
	is a hypermatrix of order $r+s$.
	
	Construct
	$
	M_{\Omega}^{{\bf j_s}\times {\bf i_r}}
	$,
	where
	$
	{\bf j_s}:=(j_1,\cdots,j_s)$, ${\bf i_r}:=(i_1,\cdots,i_r).
	$

	For $\omega_1,\cdots,\omega_s\in V^*$, $x_1,\cdots,x_r\in V$, we have	
	$$
	T(x_1,\cdots,x_r;\omega_1,\cdots,\omega_s)=\ltimes_{j=1}^q \omega_{s+1-j} M^{{\bf j_s}\times {\bf i_r}}_{\Omega}\ltimes_{i=1}^p x_i.
	$$
\end{exa}

%As an application, we consider the following example.

%%%%%%%%%%%%%%%%%%%%%%%%%%%%%%%%%%%%%%%%%%%%%%%%%%%%%%%%%%%%%%%%%%%%%%%%%%%%%%%%%%%%%%%%%%%%%%%%%%%%%%%%

\begin{exa}\label{e2.3.6}\cite{lim13} In statistical mechanics, the Yang-Baxter equation is given as follows: Let $R=[r_{i_1,i_2,i_3,i_4}]\in \F^{N\times N\times N\times N}$. Then we have (in our notation):
	\begin{align}\label{2.3.16}%{2.3.10}
		\left(R\times^4_1 R\right)\times^{(2,4)}_{(3,4)}R=R\times^{(1,2)}_{(3,4)}\times \left(R\times^4_1 R\right).
	\end{align}
	
	We express (\ref{2.3.16}) into matrix form as follows:
	
	Expressing $R\times^4_1 R$ into matrix form, we have
	\begin{align}\label{2.3.17}
		M_T:=M_R^{(1,2,3)\times 4}M_R^{1\times (2,3,4)}\in \F^{N^3\times N^3},
	\end{align}
	where its rows are indexed by $\ID(i_1,i_2,i_3;N,N,N)$, which represents the original $(i_1,i_2,i_3)$ of $R$; its columns are indexed by $\ID(i_4,i_5,i_6;N,N,N)$, which represents the original $(i_2,i_3,i_4)$ of $R$.
	
	For the product on the LHS of (\ref{2.3.16}), we convert $M_T$ into $M_T^{(i_1,i_3,i_4,i_5)\times (i_2,i_6)}$. Define
	$$
	\sigma_s:\langle 6\rangle\ra (1,3,4,5,2,6).
	$$
	Then we have
	$$
	M_T^{(1,3,4,5)\times (2,6)}={\bf V}_r^{N^2}\left[ \left( M_T \d_{N^3}^I\right)^T W^{\sigma^{-1}_s}_{\underbrace{[N\cdots N]}_6}\right].
	$$
	Hence, the LHS of (\ref{2.3.16}) becomes
	\begin{align}\label{2.3.18}
		LHS=\left({\bf V}_r^{N^2}\left[ \left( M_T \d_{N^3}^I\right)^T W^{\sigma^{-1}_s}_{\underbrace{[N\cdots N]}_6}\right]\right)M_R^{(34)\times (12)}.
	\end{align}
	
	Similarly, for the product on the LHS of (\ref{2.3.16}), we convert $M_T$ into $M_T^{(i_3,i_4)\times (i_1,i_2,i_5,i_6)}$. Define
	$$
	\sigma_t:\langle 6\rangle\ra (3,4,1,2,5,6).
	$$
	Then we have
	$$
	M_T^{(3,4)\times (1,2,5,6)}={\bf V}_r^{N^4}\left[ \left( M_T \d_{N^3}^I\right)^T W^{\sigma^{-1}_t}_{\underbrace{[N\cdots N]}_6}\right].
	$$
	Hence, the RHS of (\ref{2.3.16}) becomes
	\begin{align}\label{2.3.19}
		RHS=M_R^{(34)\times (12)}
		\left({\bf V}_r^{N^4}\left[ \left( M_T \d_{N^3}^I\right)^T W^{\sigma^{-1}_t}_{\underbrace{[N\cdots N]}_6}\right]\right).
	\end{align}
	Hence, (\ref{2.3.16}) implies that
	$$
	\begin{array}{l}
		\left({\bf V}_r^{N^2}\left[ \left( M_T \d_{N^3}^I\right)^T W^{\sigma^{-1}_s}_{\underbrace{[N\cdots N]}_6}\right]\right)M_R^{(34)\times (12)}\\
		=M_R^{(34)\times (12)}\left({\bf V}_r^{N^4}\left[ \left( M_T \d_{N^3}^I\right)^T W^{\sigma^{-1}_t}_{\underbrace{[N\cdots N]}_6}\right]\right).
	\end{array}
	$$
\end{exa}

\section{Conclusion}

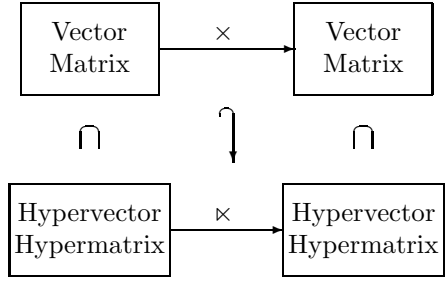
\begin{figure}
	\centering
	\setlength{\unitlength}{6 mm}
	\begin{picture}(10,7)(-0.5,-0.5)
		%\put(0,0){\framebox(3,2){$\begin{array}{c}\mbox{\begin{tiny}Hypervector\end{tiny}}\\
					%\mbox{\begin{tiny}Hypermatrix\end{tiny}}\end{array}$}}
		%\put(6,0){\framebox(3,2){$\begin{array}{c}\mbox{\begin{tiny}Hypervector\end{tiny}}\\
					%\mbox{\begin{tiny}Hypermatrix\end{tiny}}\end{array}$}}
		\put(-0.25,0){\framebox(3.5,2){$\begin{array}{c}\mbox{Hypervector}\\
					\mbox{Hypermatrix}\end{array}$}}
		\put(5.75,0){\framebox(3.5,2){$\begin{array}{c}\mbox{Hypervector}\\
					\mbox{Hypermatrix}\end{array}$}}
		\put(0,4){\framebox(3,2){$\begin{array}{c}\mbox{Vector}\\
					\mbox{Matrix}\end{array}$}}
		\put(6,4){\framebox(3,2){$\begin{array}{c}\mbox{Vector}\\
					\mbox{Matrix}\end{array}$}}
		\put(3.25,1){\vector(1,0){2.5}}
		\put(3,5){\vector(1,0){3}}
		\put(4.2,1.2){$\ltimes$}
		\put(4.2,5.2){$\times$}
		\put(1.5,2.8){\oval(0.4,1)[t]}
		\put(7.5,2.8){\oval(0.4,1)[t]}
		\put(4.5,3.5){\oval(0.3,0.4)[t]}
		\put(4.65,3.5){\vector(0,-1){1}}
	\end{picture}
	\caption{$\times$ vs $\ltimes$ \label{Fig.c.1}}
\end{figure}

In this paper, the matrix expression of hypermatrices was first proposed. As an auxiliary tool, the permutation matrices have also been discussed in detail. It is used to reveal certain properties of the matrix expression of hypermatrices. Then it was shown that the STPs of matrices are essentially the multilinear operators of hypermatrices (including hypervectors). The operators over hypermatrices including contracted products, are realized by STPs through the matrix expression of hypermatrices. In fact, the actions of STP over hypermatrices can be considered as a generalization of the actions of the conventional matrix product over matrices. This fact can be demonstrated by Figure \ref{Fig.c.1}.

\section*{Appendix}

In the following some permutation matrices are listed.

\begin{onecolumn}

\begin{enumerate}
\item $d=3$, $n=2$:

\begin{itemize}
\item[(i)] $\sigma_1:[1,2,3]\ra [1,2,3]$:
$$
W_2^{\sigma_1}=\d_8[1,2,3,4,5,6,7,8].
$$
\item[(ii)] $\sigma_2:[1,2,3]\ra [1,3,2]$:
$$
W_2^{\sigma_2}=\d_8[1,3,2,4,5,7,6,8].
$$
\item[(iii)] $\sigma_3:[1,2,3]\ra [2,1,3]$:
$$
W_2^{\sigma_3}=\d_8[1,2,5,6,3,4,7,8].
$$
\item[(iv)] $\sigma_4:[1,2,3]\ra [2,3,1]$:
$$
W_2^{\sigma_4}=\d_8[1,3,5,7,2,4,6,8]..
$$
\item[(v)] $\sigma_5:[1,2,3]\ra [3,1,2]$:
$$
W_2^{\sigma_5}=\d_8[1,5,2,6,3,7,4,8].
$$
\item[(vi)] $\sigma_6:[1,2,3]\ra [3,2,1]$:
$$
W_2^{\sigma_6}=\d_8[1,5,3,7,2,6,4,8].
$$
\end{itemize}

\item $d=3$, $n=3$:

\begin{itemize}

\item[(i)] $\sigma_1:[1,2,3]\ra [1,2,3]$:
$$
W_3^{\sigma_1}=\d_{27}[ 1, 2, 3, 4, 5, 6, 7, 8, 9,10,11,12,13,14,15,16,17,18,19,20,21,22,23,24,25,26,27].
$$
\item[(ii)] $\sigma_2:[1,2,3]\ra [1,3,2]$:
$$
W_3^{\sigma_2}=\d_{27}[1,4,7,2,5,8,3,6,9,10,13,16,11,14,17,12,15,18,19,22,25,20,23,26,21,24,27].
$$
\item[(iii)] $\sigma_3:[1,2,3]\ra [2,1,3]$:
$$
W_3^{\sigma_3}=\d_{27}[1,2,3,10,11,12,19,20,21,4,5,6,13,14,15,22,23,24,7,8,9,16,17,18,25,26,27].
$$
\item[(iv)] $\sigma_4:[1,2,3]\ra [2,3,1]$:
$$
W_3^{\sigma_4}=\d_{27}[1,10,19,2,11,20,3,12,21,4,13,22,5,14,23,6,15,24,7,16,25,8,17,26,9,18,27].
$$

\item[(v)] $\sigma_5:[1,2,3]\ra [3,1,2]$:
$$
W_3^{\sigma_5}=\d_{27}[1,4,7,10,13,16,19,22,25,2,5,8,11,14,17,20,23,26,3,6,9,12,15,18,21,24,27].
$$
\item[(vi)] $\sigma_6:[1,2,3]\ra [3,2,1]$:
$$
W_3^{\sigma_6}=\d_{27}[1,10,19,4,13,22,7,16,25,2,11,20,5,14,23,8,17,26,3,12,21,6,15,24,9,18,27].
$$
\end{itemize}

\item $d=3$, $n=4$:
\begin{itemize}
\item[(i)] $\sigma_1:[1,2,3]\ra [1,2,3]$:
$$
\begin{array}{l}
W_4^{\sigma_1}=\d_{64}[1,2,3,4,5,6,7,8,9,10,11,12,13,14,15,16,17,18,19,20,21,22,\\
~~~~23,24,25,26,27,28,29,30,31,32,33,34,35,36,37,38,39,40,41,42,43,\\
~~~~44,45,46,47,48,49,50,51,52,53,54,55,56,57,58,59,60,61,62,63,64].
\end{array}
$$
\item[(ii)] $\sigma_2:[1,2,3]\ra [1,3,2]$:
$$
\begin{array}{l}
W_4^{\sigma_2}=\d_{64}[1,5,9,13,2,6,10,14,3,7,11,15,4,8,12,16,17,21,25, 29,18,22,\\
~~~~26,30,19,23,27,31,20,24,28,32,33,37,41,45,34,38,42,46,35,39,43,\\
~~~~47,36,40,44,48,49,53,57,61,50,54,58,62,51,55,59,63,52,56,60,64].
\end{array}
$$
\item[(iii)] $\sigma_3:[1,2,3]\ra [2,1,3]$:
$$
\begin{array}{l}
W_4^{\sigma_3}=\d_{64}[1,2,3,4,17,18,19,20,33,34,35,36,49,50,51,52,5,6,7,8,21,22,\\
~~~~23,24,37,38,39,40,53,54,55,56,9,10,11,12,25,26,27,28,41,42,43,44,\\
~~~~57,58,59,60,13,14,15,16,29,30,31,32,45,46,47,48,61,62,63,64].
\end{array}
$$
\item[(iv)] $\sigma_4:[1,2,3]\ra [2,3,1]$:
$$
\begin{array}{l}
W_4^{\sigma_4}=\d_{64}[1,17,33,49,2,18,34,50,3,19,35,51,4,20,36,52,5,21,37,53,6,\\
~~~~22,38,54,7,23,39,55,8,24,40,56,9,25,41,57,10,26,42,58,11,27,43,\\
~~~~59,12,28,44,60,13,29,45,61,14,30,46,62,15,31,47,63,16,32,48,64].
\end{array}
$$
\item[(v)] $\sigma_5:[1,2,3]\ra [3,1,2]$:
$$
\begin{array}{l}
W_5^{\sigma_4}=\d_{64}[1,5,9,13,17,21,25,29,33,37,41,45,49,53,57,61,2,6,10,14,\\
~~~~18,22,26,30,34,38,42,46,50,54,58,62,3,7,11,15,19,23,27,31,35,39,\\
~~~~43,47,51,55,59,63,4,8,12,16,20,24,28,32,36,40,44,48,52,56,60,64].
\end{array}
$$
\item[(vi)] $\sigma_6:[1,2,3]\ra [3,2,1]$:
$$
\begin{array}{l}
W_6^{\sigma_4}=\d_{64}[1,17,33,49,5,21,37,53,9,25,41,57,13,29,45,61,2,18,34,50,\\
~~~~6,22,38,54,10,26,42,58,14,30,46,62,3,19,35,51,7,23,39,55,11,27,\\
~~~~43,59,15,31,47,63,4,20,36,52,8,24,40,56,12,28,44,60,16,32,48,64].
\end{array}
$$
\end{itemize}
\item $d=4$, $n=2$:
\begin{itemize}
\item[(i)] $\sigma_1:[1,2,3,4]\ra [1,2,3,4]$:
$$
W_2^{\sigma_1}=\d_{16}[1,2,3,4,5,6,7,8,9,10,11,12,13,14,15,16].
$$
\item[(ii)] $\sigma_2:[1,2,3,4]\ra [1,2,4,3]$:
$$
W_2^{\sigma_2}=\d_{16}[1,3,2,4,5,7,6,8,9,11,10,12,13,15,14,16].
$$
\item[(iii)] $\sigma_3:[1,2,3,4]\ra [1,3,2,4]$:
$$
W_2^{\sigma_3}=\d_{16}[1,2,5,6,3,4,7,8,9,10,13,14,11,12,15,16].
$$
\item[(iv)] $\sigma_4:[1,2,3,4]\ra [1,3,4,2]$:
$$
W_2^{\sigma_4}=\d_{16}[1,5,2,6,3,7,4,8,9,13,10,14,11,15,12,16].
$$
\item[(v)] $\sigma_5:[1,2,3,4]\ra [1,4,2,3]$:
$$
W_2^{\sigma_5}=\d_{16}[1,3,5,7,2,4,6,8,9,11,13,15,10,12,14,16].
$$
\item[(vi)] $\sigma_6:[1,2,3,4]\ra [1,4,3,2]$:
$$
W_2^{\sigma_6}=\d_{16}[1,5,3,7,2,6,4,8,9,13,11,15,10,14,12,16].
$$
\item[(vii)] $\sigma_7:[1,2,3,4]\ra [2,1,3,4]$:
$$
W_2^{\sigma_7}=\d_{16}[1,2,3,4,9,10,11,12,5,6,7,8,13,14,15,16].
$$
\item[(viii)] $\sigma_8:[1,2,3,4]\ra [2,1,4,3]$:
$$
W_2^{\sigma_8}=\d_{16}[1,3,2,4,9,11,10,12,5,7,6,8,13,15,14,16].
$$
\item[(ix)] $\sigma_9:[1,2,3,4]\ra [2,3,1,4]$:
$$
W_2^{\sigma_8}=\d_{16}[1,2,9,10,3,4,11,12,5,6,13,14,7,8,15,16].
$$
\item[(x)] $\sigma_{10}:[1,2,3,4]\ra [2,3,4,1]$:
$$
W_2^{\sigma_{10}}=\d_{16}[1,9,2,10,3,11,4,12,5,13,6,14,7,15,8,16].
$$
\item[(xi)] $\sigma_{11}:[1,2,3,4]\ra [2,4,1,3]$:
$$
W_2^{\sigma_{11}}=\d_{16}[1,3,9,11,2,4,10,12,5,7,13,15,6,8,14,16].
$$
\item[(xii)] $\sigma_{12}:[1,2,3,4]\ra [2,4,3,1]$:
$$
W_2^{\sigma_{12}}=\d_{16}[1,9,3,11,2,10,4,12,5,13,7,15,6,14,8,16].
$$
\item[(xiii)] $\sigma_{13}:[1,2,3,4]\ra [3,1,2,4]$:
$$
W_2^{\sigma_{13}}=\d_{16}[1,2,5,6,9,10,13,14,3,4,7,8,11,12,15,16].
$$
\item[(xiv)] $\sigma_{14}:[1,2,3,4]\ra [3,1,4,2]$:
$$
W_2^{\sigma_{14}}=\d_{16}[1,5,2,6,9,13,10,14,3,7,4,8,11,15,12,16].
$$
\item[(xv)] $\sigma_{15}:[1,2,3,4]\ra [3,2,1,4]$:
$$
W_2^{\sigma_{15}}=\d_{16}[1,5,2,6,9,13,10,14,3,7,4,8,11,15,12,16].
$$
\item[(xvi)] $\sigma_{16}:[1,2,3,4]\ra [3,2,4,1]$:
$$
W_2^{\sigma_{16}}=\d_{16}[1,9,2,10,5,13,6,14,3,11,4,12,7,15,8,16].
$$
\item[(xvii)] $\sigma_{17}:[1,2,3,4]\ra [3,4,1,2]$:
$$
W_2^{\sigma_{17}}=\d_{16}[1,5,9,13,2,6,10,14,3,7,11,15,4,8,12,16].
$$
\item[(xviii)] $\sigma_{18}:[1,2,3,4]\ra [3,4,2,1]$:
$$
W_2^{\sigma_{18}}=\d_{16}[1,9,5,13,2,10,6,14,3,11,7,15,4,12,8,16].
$$
\item[(xix)] $\sigma_{19}:[1,2,3,4]\ra [4,1,2,3]$:
$$
W_2^{\sigma_{19}}=\d_{16}[1,3,5,7,9,11,13,15,2,4,6,8,10,12,14,16].
$$
\item[(xx)] $\sigma_{20}:[1,2,3,4]\ra [4,1,3,2]$:
$$
W_2^{\sigma_{20}}=\d_{16}[1,5,3,7,9,13,11,15,2,6,4,8,10,14,12,16].
$$
\item[(xxi)] $\sigma_{21}:[1,2,3,4]\ra [4,2,1,3]$:
$$
W_2^{\sigma_{21}}=\d_{16}[1,3,9,11,5,7,13,15,2,4,10,12,6,8,14,16].
$$
\item[(xxii)] $\sigma_{22}:[1,2,3,4]\ra [4,2,3,1]$:
$$
W_2^{\sigma_{22}}=\d_{16}[1,9,3,11,5,13,7,15,2,10,4,12,6,14,8,16].
$$
\item[(xxiii)] $\sigma_{23}:[1,2,3,4]\ra [4,3,1,2]$:
$$
W_2^{\sigma_{23}}=\d_{16}[1,5,9,13,3,7,11,15,2,6,10,14,4,8,12,16].
$$
\item[(xxiv)] $\sigma_{24}:[1,2,3,4]\ra [4,3,2,1]$:
$$
W_2^{\sigma_{24}}=\d_{16}[1,9,5,13,3,11,7,15,2,10,6,14,4,12,8,16].
$$
\end{itemize}

\item $d=4$, $n=3$:
\begin{itemize}
\item[(i)] $\sigma_1:[1,2,3,4]\ra [1,2,3,4]$:
$$
\begin{array}{l}
W_3^{\sigma_1}=\d_{81}[1,2,3,4,5,6,7,8,9,10,11,12,13,14,15,16,17,18,19,20,\\
~~21,22,23,24,25,26,27,28,29,30,31,32,33,34,35,36,37,38,39,40,\\
~~41,42,43,44,45,46,47,48,49,50,51,52,53,54,55,56,57,58,59,60,\\
~~61,62,63,64,65,66,67,68,69,70,71,72,73,74,75,76,77,78,79,80,81].
\end{array}
$$
\item[(ii)] $\sigma_2:[1,2,3,4]\ra [1,2,4,3]$:
$$
\begin{array}{l}
W_3^{\sigma_2}=\d_{81}[1,4,7,2,5,8,3,6,9,10,13,16,11,14,17,12,15,18,19,22,\\
~~25,20,23,26,21,24,27,28,31,34,29,32,35,30,33,36,37,40,43,38,\\
~~41,44,39,42,45,46,49,52,47,50,53,48,51,54,55,58,61,56,59,62,\\
~~57,60,63,64,67,70,65,68,71,66,69,72,73,76,79,74,77,80,75,78,81].
\end{array}
$$
\item[(iii)] $\sigma_3:[1,2,3,4]\ra [1,3,2,4]$:
$$
\begin{array}{l}
W_3^{\sigma_3}=\d_{81}[1,2,3,10,11,12,19,20,21,4,5,6,13,14,15,22,23,24,7,8,\\
~~9,16,17,18,25,26,27,28,29,30,37,38,39,46,47,48,31,32,33,40,\\
~~41,42,49,50,51,34,35,36,43,44,45,52,53,54,55,56,57,64,65,66,\\
~~73,74,75,58,59,60,67,68,69,76,77,78,61,62,63,70,71,72,79,80,81].
\end{array}
$$
\item[(iv)] $\sigma_4:[1,2,3,4]\ra [1,3,4,2]$:
$$
\begin{array}{l}
W_3^{\sigma_4}=\d_{81}[1,10,19,2,11,20,3,12,21,4,13,22,5,14,23,6,15,24,7,16,\\
~~25,8,17,26,9,18,27,28,37,46,29,38,47,30,39,48,31,40,49,32,\\
~~41,50,33,42,51,34,43,52,35,44,53,36,45,54,55,64,73,56,65,74,\\
~~57,66,75,58,67,76,59,68,77,60,69,78,61,70,79,62,71,80,63,72,81].
\end{array}
$$
\item[(v)] $\sigma_5:[1,2,3,4]\ra [1,4,2,3]$:
$$
\begin{array}{l}
W_3^{\sigma_5}=\d_{81}[1,4,7,10,13,16,19,22,25,2,5,8,11,14,17,20,23,26,3,6,\\
~~9,12,15,18,21,24,27,28,31,34,37,40,43,46,49,52,29,32,35,38,\\
~~41,44,47,50,53,30,33,36,39,42,45,48,51,54,55,58,61,64,67,70,\\
~~73,76,79,56,59,62,65,68,71,74,77,80,57,60,63,66,69,72,75,78,81].
\end{array}
$$
\item[(vi)] $\sigma_6:[1,2,3,4]\ra [1,4,3,2]$:
$$
\begin{array}{l}
W_3^{\sigma_6}=\d_{81}[1,10,19,4,13,22,7,16,25,2,11,20,5,14,23,8,17,26,3,12,\\
~~21,6,15,24,9,18,27,28,37,46,31,40,49,34,43,52,29,38,47,32,\\
~~41,50,35,44,53,30,39,48,33,42,51,36,45,54,55,64,73,58,67,76,\\
~~61,70,79,56,65,74,59,68,77,62,71,80,57,66,75,60,69,78,63,72,81].
\end{array}
$$
\item[(vii)] $\sigma_7:[1,2,3,4]\ra [2,1,3,4]$:
$$
\begin{array}{l}
W_3^{\sigma_7}=\d_{81}[1,2,3,4,5,6,7,8,9,28,29,30,31,32,33,34,35,36,55,56,\\
~~57,58,59,60,61,62,63,10,11,12,13,14,15,16,17,18,37,38,39,40,\\
~~41,42,43,44,45,64,65,66,67,68,69,70,71,72,19,20,21,22,23,24,\\
~~25,26,27,46,47,48,49,50,51,52,53,54,73,74,75,76,77,78,79,80,81].
\end{array}
$$
\item[(viii)] $\sigma_8:[1,2,3,4]\ra [2,1,4,3]$:
$$
\begin{array}{l}
W_3^{\sigma_8}=\d_{81}[1,4,7,2,5,8,3,6,9,28,31,34,29,32,35,30,33,36,55,58,\\
~~61,56,59,62,57,60,63,10,13,16,11,14,17,12,15,18,37,40,43,38,\\
~~41,44,39,42,45,64,67,70,65,68,71,66,69,72,19,22,25,20,23,26,\\
~~21,24,27,46,49,52,47,50,53,48,51,54,73,76,79,74,77,80,75,78,81].
\end{array}
$$

\item[(ix)] $\sigma_9:[1,2,3,4]\ra [2,3,1,4]$:
$$
\begin{array}{l}
W_3^{\sigma_8}=\d_{81}[1,2,3,28,29,30,55,56,57,4,5,6,31,32,33,58,59,60,7,8,\\
~~9,34,35,36,61,62,63,10,11,12,37,38,39,64,65,66,13,14,15,40,\\
~~41,42,67,68,69,16,17,18,43,44,45,70,71,72,19,20,21,46,47,48,\\
~~73,74,75,22,23,24,49,50,51,76,77,78,25,26,27,52,53,54,79,80,81].
\end{array}
$$

\item[(x)] $\sigma_{10}:[1,2,3,4]\ra [2,3,4,1]$:
$$
\begin{array}{l}
W_3^{\sigma_{10}}=\d_{81}[1,28,55,2,29,56,3,30,57,4,31,58,5,32,59,6,33,60,7,34,\\
~~61,8,35,62,9,36,63,10,37,64,11,38,65,12,39,66,13,40,67,14,\\
~~41,68,15,42,69,16,43,70,17,44,71,18,45,72,19,46,73,20,47,74,\\
~~21,48,75,22,49,76,23,50,77,24,51,78,25,52,79,26,53,80,27,54,81].
\end{array}
$$
\item[(xi)] $\sigma_{11}:[1,2,3,4]\ra [2,4,1,3]$:
$$
\begin{array}{l}
W_3^{\sigma_{11}}=\d_{81}[1,4,7,28,31,34,55,58,61,2,5,8,29,32,35,56,59,62,3,6,\\
~~9,30,33,36,57,60,63,10,13,16,37,40,43,64,67,70,11,14,17,38,\\
~~41,44,65,68,71,12,15,18,39,42,45,66,69,72,19,22,25,46,49,52,\\
~~73,76,79,20,23,26,47,50,53,74,77,80,21,24,27,48,51,54,75,78,81].
\end{array}
$$
\item[(xii)] $\sigma_{12}:[1,2,3,4]\ra [2,4,3,1]$:
$$
\begin{array}{l}
W_3^{\sigma_{12}}=\d_{81}[1,28,55,4,31,58,7,34,61,2,29,56,5,32,59,8,35,62,3,30,\\
~~57,6,33,60,9,36,63,10,37,64,13,40,67,16,43,70,11,38,65,14,\\
~~41,68,17,44,71,12,39,66,15,42,69,18,45,72,19,46,73,22,49,76,\\
~~25,52,79,20,47,74,23,50,77,26,53,80,21,48,75,24,51,78,27,54,81].
\end{array}
$$

\item[(xiii)] $\sigma_{13}:[1,2,3,4]\ra [3,1,2,4]$:
$$
\begin{array}{l}
W_3^{\sigma_{13}}=\d_{81}[1,2,3,10,11,12,19,20,21,28,29,30,37,38,39,46,47,48,55,56,\\
~~57,64,65,66,73,74,75,4,5,6,13,14,15,22,23,24,31,32,56,57,\\
~~64,65,66,73,74,75,4,5,6,13,14,15,22,23,24,31,32,16,17,18,\\
~~25,26,27,34,35,36,43,44,45,52,53,54,61,62,63,70,71,72,79,80,81].
\end{array}
$$

\item[(xiv)] $\sigma_{14}:[1,2,3,4]\ra [3,1,4,2]$:
$$
\begin{array}{l}
W_3^{\sigma_{14}}=\d_{81}[1,10,19,2,11,20,3,12,21,28,37,46,29,38,47,30,39,48,55,64,\\
~~73,56,65,74,57,66,75,4,13,22,5,14,23,6,15,24,31,40,49,32,\\
~~41,50,33,42,51,58,67,76,59,68,77,60,69,78,7,16,25,8,17,26,\\
~~9,18,27,34,43,52,35,44,53,36,45,54,61,70,79,62,71,80,63,72,81].
\end{array}
$$

\item[(xv)] $\sigma_{15}:[1,2,3,4]\ra [3,2,1,4]$:
$$
\begin{array}{l}
W_3^{\sigma_{15}}=\d_{81}[1,10,19,2,11,20,3,12,21,28,37,46,29,38,47,30,39,48,55,64,\\
~~73,56,65,74,57,66,75,4,13,22,5,14,23,6,15,24,31,40,49,32,\\
~~41,50,33,42,51,58,67,76,59,68,77,60,69,78,7,16,25,8,17,26,\\
~~9,18,27,34,43,52,35,44,53,36,45,54,61,70,79,62,71,80,63,72,81].
\end{array}
$$

\item[(xvi)] $\sigma_{16}:[1,2,3,4]\ra [3,2,4,1]$:
$$
\begin{array}{l}
W_3^{\sigma_{16}}=\d_{16}[1,28,55,2,29,56,3,30,57,10,37,64,11,38,65,12,39,66,19,46,\\
~~73,20,47,74,21,48,75,4,31,58,5,32,59,6,33,60,13,40,67,14,\\
~~41,68,15,42,69,22,49,76,23,50,77,24,51,78,7,34,61,8,35,62,\\
~~9,36,63,16,43,70,17,44,71,18,45,72,25,52,79,26,53,80,27,54,81].
\end{array}
$$

\item[(xvii)] $\sigma_{17}:[1,2,3,4]\ra [3,4,1,2]$:
$$
\begin{array}{l}
W_3^{\sigma_{17}}=\d_{81}[1,10,19,28,37,46,55,64,73,2,11,20,29,38,47,56,65,74,3,12,\\
~~21,30,39,48,57,66,75,4,13,22,31,40,49,58,67,76,5,14,23,32,\\
~~41,50,59,68,77,6,15,24,33,42,51,60,69,78,7,16,25,34,43,52,\\
~~61,70,79,8,17,26,35,44,53,62,71,80,9,18,27,36,45,54,63,72,81].
\end{array}
$$

\item[(xviii)] $\sigma_{18}:[1,2,3,4]\ra [3,4,2,1]$:
$$
\begin{array}{l}
W_3^{\sigma_{18}}=\d_{81}[1,28,55,10,37,64,19,46,73,2,29,56,11,38,65,20,47,74,3,30,\\
~~57,12,39,66,21,48,75,4,31,58,13,40,67,22,49,76,5,32,59,14,\\
~~41,68,23,50,77,6,33,60,15,42,69,24,51,78,7,34,61,16,43,70,\\
~~25,52,79,8,35,62,17,44,71,26,53,80,9,36,63,18,45,72,27,54,81].
\end{array}
$$

\item[(xix)] $\sigma_{19}:[1,2,3,4]\ra [4,1,2,3]$:
$$
\begin{array}{l}
W_3^{\sigma_{19}}=\d_{81}[1,4,7,10,13,16,19,22,25,28,31,34,37,40,43,46,49,52,55,58,\\
~~61,64,67,70,73,76,79,2,5,8,11,14,17,20,23,26,29,32,35,38,\\
~~41,44,47,50,53,56,59,62,65,68,71,74,77,80,3,6,9,12,15,18,\\
~~21,24,27,30,33,36,39,42,45,48,51,54,57,60,63,66,69,72,75,78,81].
\end{array}
$$

\item[(xx)] $\sigma_{20}:[1,2,3,4]\ra [4,1,3,2]$:
$$
\begin{array}{l}
W_3^{\sigma_{20}}=\d_{81}[1,10,19,4,13,22,7,16,25,28,37,46,31,40,49,34,43,52,55,64,\\
~~73,58,67,76,61,70,79,2,11,20,5,14,23,8,17,26,29,38,47,32,\\
~~41,50,35,44,53,56,65,74,59,68,77,62,71,80,3,12,21,6,15,24,\\
~~9,18,27,30,39,48,33,42,51,36,45,54,57,66,75,60,69,78,63,72,81].
\end{array}
$$

\item[(xxi)] $\sigma_{21}:[1,2,3,4]\ra [4,2,1,3]$:
$$
\begin{array}{l}
W_3^{\sigma_{21}}=\d_{81}[1,4,7,28,31,34,55,58,61,10,13,16,37,40,43,64,67,70,19,22,\\
~~25,46,49,52,73,76,79,2,5,8,29,32,35,56,59,62,11,14,17,38,\\
~~41,44,65,68,71,20,23,26,47,50,53,74,77,80,3,6,9,30,33,36,\\
~~57,60,63,12,15,18,39,42,45,66,69,72,21,24,27,48,51,54,75,78,81].
\end{array}
$$

\item[(xxii)] $\sigma_{22}:[1,2,3,4]\ra [4,2,3,1]$:
$$
\begin{array}{l}
W_3^{\sigma_{22}}=\d_{81}[1,28,55,4,31,58,7,34,61,10,37,64,13,40,67,16,43,70,19,46,\\
~~73,22,49,76,25,52,79,2,29,56,5,32,59,8,35,62,11,38,65,14,\\
~~41,68,17,44,71,20,47,74,23,50,77,26,53,80,3,30,57,6,33,60,\\
~~9,36,63,12,39,66,15,42,69,18,45,72,21,48,75,24,51,78,27,54,81].
\end{array}
$$

\item[(xxiii)] $\sigma_{23}:[1,2,3,4]\ra [4,3,1,2]$:
$$
\begin{array}{l}
W_3^{\sigma_{23}}=\d_{81}[1,10,19,28,37,46,55,64,73,4,13,22,31,40,49,58,67,76,7,16,\\
~~25,34,43,52,61,70,79,2,11,20,29,38,47,56,65,74,5,14,23,32,\\
~~41,50,59,68,77,8,17,26,35,44,53,62,71,80,3,12,21,30,39,48,\\
~~57,66,75,6,15,24,33,42,51,60,69,78,9,18,27,36,45,54,63,72,81].
\end{array}
$$

\item[(xxiv)] $\sigma_{24}:[1,2,3,4]\ra [4,3,2,1]$:
$$
\begin{array}{l}
W_3^{\sigma_{24}}=\d_{81}[1,28,55,10,37,64,19,46,73,4,31,58,13,40,67,22,49,76,7,34,\\
~~61,16,43,70,25,52,79,2,29,56,11,38,65,20,47,74,5,32,59,14,\\
~~41,68,23,50,77,8,35,62,17,44,71,26,53,80,3,30,57,12,39,66,\\
~~21,48,75,6,33,60,15,42,69,24,51,78,9,36,63,18,45,72,27,54,81].
\end{array}
$$
\end{itemize}
\end{enumerate}
\end{onecolumn}


\begin{thebibliography}{00}
		%\bibitem{bar22}  E. Bar-Shalom, O. Dalin, M. Margaliot, Compound matrices in systems and control theory: a tutorial, {\it arXiv:2204.00676v1},  2022.
		%
		%\bibitem{boo86} W.M. Boothby, {\it Introduction to Differentiable Manifolds and Riemannian Geometry}, 2nd Ed., Elsevier, London, 1986.
		
		\bibitem{bus19} B. Bush, J. Culp, K. Pearson, {errpm-Frobenius theorem for hypermatrices in the max algebra, {\it Discrete Math.}, Vol. 342, 64-73, 2019.
			%
			\bibitem{che01} D. Cheng, Semi-tensor product of matrices and its application to Morgen's problem, {\it Science in China}, Series F, Vol. 44, No. 3, 195-212,2001.
			%
			%\bibitem{che07} D. Cheng, H. Qi, {\it Semi-tensor Product of Matrices - Theory and Applications}, Science Press, Beijing, 2007.
			%
            \bibitem{che07} D. Cheng, H. Qi, {\it Semi-tensor Product of Matrices - Theory and Applications}, Science Press, Beijing, 2007, (in Chinese)
   
			\bibitem{che11} D. Cheng, H. Qi, Z. Li, {\it Analysis and Control of Boolean Networks - A Semi-tensor Product
				Approach}, Springer, London, 2011.
			%%
			\bibitem{che12}  D. Cheng, H. Qi, Y. Zhao, {\it An Introduction to Semi-tensor Product of Matrices and Its Applications}, World Scientific, Singapore, 2012.
			%
			\bibitem{che15}  D. Cheng, F. He, H. Qi, T. Xu, Modeling, analysis and control of networked evolutionary games, {\it IEEE Trans. Aut. Contr.}, Vol. 60, No. 9, 2402-2415, 2015.
			%
			\bibitem{che19} D. Cheng, On equivalence of matrices, {\it Asian Journal of Math.}, Vol. 23, No. 2, 257-348, 2019.
			%
			\bibitem{che19b} D. Cheng, {\it From Dimension-Free Matrix Theory to Cross-Dimensional Dynamic Systems}, Elsevier, London, 2019.
			%
			\bibitem{che19c} D. Cheng, Z. Liu, A new semi-tensor product of matrices, {\it Contr. Theory Tech.}, Vol. 17, No. 1, 14-22, 2019.
			
			\bibitem{che20} D. Cheng, H. Qi, {\it Lecture Notes in Semi-Tensor Product of Matrices}, Vol. 1, {\it Basic Theory and Multilinear Operation}, Science Press, Beijing, 2020, (in Chinese).
			%
			\bibitem{che21} D. Cheng, Y. Wu, G. Zhao, S. Fu, A comprehensive survey on STP approach to finite games, {\it J. Sys. Sci. Compl.}, Vol. 34, No. 5, 1666-1680, 2021.
			
			\bibitem{che21b}  D. Cheng, Z. Ji, J. Feng, S. Fu, J. Zhao, Perfect hypercomplex algebras: Semi-tensor product approach, {\it Math. Model. Contr.}, Vol. 1, No. 4, 177-187, 2021.
			%
			\bibitem{chepr} D. Cheng, Z. Ji, From dimension-free manifold to dimension-varying control system,  {\it Commun. Inform. Sys.}, Vol. 23, No. 1, 85-150, 2023.
			%
			\bibitem{che23} D. Cheng,  X. Zhang, Z. Ji, Semi-tensor product of hypermatrices with application to compound hypermatrices,  http:arxiv.org/abs/2303.06295, 2023.
			%
			\bibitem{che23b} D. Cheng, Z. Ji,  {\it Lecture Notes in Semi-Tensor Product of Matrices}, Vol. 4, {\it Finite and Dimensional-free Dynamic Systems}, Science Press, Beijing, 2023, (in Chinese).
			%
			\bibitem{del00} L. De Lathauwer, B. De Moor, J. Vandewalle, On the best rak-1 and rank-(R1,R2,$\cdots$,Rn) approximation of higher order tensors, {\it SIAM J. Matrix Anal. Appl.}, Vol. 21, 1324-1342, 2000.
			%
			\bibitem{del01} L. De Lathauwer, B. De Moor, J. Vandewalle, Independent component analysis and (simultaneous) third-order tensor diagonalization}, {\it IEEE Trans. Signal Proc.}, Vol. 49, 2262-2271, 2001.
		%
		\bibitem{fan20} Z. Fan, C. Deng, H. Li, C. Bu, Multiplications and eignevalues of tensors via linear maps, Vol. 68, No. 3, 606-621, 2020.
		%
		\bibitem{for16} E. Fornasini, M.E. Valcher, Recent developments in Boolean networks control, {\it J. Contr. Dec.}, Vol. 3, No. 1, 1-18, 2016.
		%
		\bibitem{li18} H. Li, G. Zhao, M. Meng, J. Feng, A survey on applications of semi-tensor product method in engineering, {\it Science China}, Vol. 61, 010202:1-010202:17, 2018.
		%
		\bibitem{lim13}
		L. Lim, Tensors and Hypermatrices, in L. Hogben (Ed.) {\it Handbook of Linear Algebra} (2nd ed.), Chapter 15,
		Chapman and Hall/CRC.https://doi.org/10.1201/b16113, 2013.
		%
		\bibitem{lu17}  J. Lu, H. Li, Y. Liu, F. Li, Survey on semi-tensor product method with its applications in logical networks and other finite-valued systems, {\it IET Contr. Theory Appl.}, Vol. 11, No. 13, 2040-2047, 2017.
		%
		\bibitem{muh16} A. Muhammad, A. Rushdi, F.A. M. Ghaleb, A tutorial exposition of semi-tensor products of matrices with a stress on their representation of Boolean function, {\it JKAU Comp. Sci.}, Vol. 5, 3-30, 2016.
		%
		%\bibitem{qi05} L. Qi, Eigenvalues of a real supersymmetric tensor, {\it J. Symbolic Comput.}, Vol. 21, 1302-1324, 2005.
		%
		%\bibitem{qi07} L. Qi, Eigenvalues and invariants of tensors, {\it J. Math. Anal. Appl.}, Vol. 325, 1363-1377, 2007.
		%
		\bibitem{vas03} M.A.O. Vasilescu, D. Terzopoulos, Multilinear subspace analysis for image ensembles, {\it Proc. Cof. CVPR03}, Madison, WI, 93-99, 2003.
		%
		%\bibitem{wu22} C. Wu, R. Pines, M. Margaliot, J.J. Slotine, Generalization of the multiplicative and additive compounds of square matrices and contraction theory in the hausdorff dimension, {\it IEEE Trans. Aut. Contr.}, Vol. 67, No. 9, 4629-4644, 2022.
		
		\bibitem{yan22} Y. Yan, D. Cheng, J. Feng, H. Li, J. Yue, Survey on applications of algebraic state space theory of logical systems to finite state machines, {\it Sci. China Inf. Sci.}, https://doi.org/10.1007/s11432-22-3538-4, 2022.
	\end{thebibliography}
\end{document}